\documentclass[11pt]{amsart}
\usepackage{lmodern}
\usepackage{amsmath, amsthm, amssymb, amsfonts}
\usepackage{mathtools}
\usepackage{stmaryrd}
\usepackage[normalem]{ulem}
\usepackage{mathrsfs}
\usepackage{verbatim} 
\usepackage{longtable}
\usepackage{mathtools}
\usepackage{enumitem}
\usepackage{bbm}

\usepackage{tikz}
\usetikzlibrary{decorations.pathmorphing}
\tikzset{snake it/.style={decorate, decoration=snake}}
\usepackage{caption}
\usepackage{tikz-cd}
\usetikzlibrary{arrows}
\usepackage{hyperref}
\usepackage{url}
\usepackage[scr=boondox]{mathalfa}

\DeclareMathAlphabet{\matheu}{U}{eur}{m}{n}\usepackage{tensor}
\usepackage{hyperref, mathrsfs}
\usepackage{setspace}
\usepackage{graphicx}
\usepackage{enumitem}
\usepackage[headings]{fullpage}
\usepackage{tikz}
\usepackage{tikz-cd}
\usepackage{wrapfig}
\usepackage{float}
\makeatletter
\newcommand*\bigcdot{\mathpalette\bigcdot@{.5}}
\newcommand*\bigcdot@[2]{\mathbin{\vcenter{\hbox{\scalebox{#2}{$\m@th#1\bullet$}}}}}

\makeatletter
\let\@wraptoccontribs\wraptoccontribs

\usepackage{kantlipsum}
\newcommand\Z{\mathbb{Z}}
\usepackage{bm}
\usepackage[cmtip,all]{xy}
 
\newtheorem*{theorem*}{Theorem}

\newtheorem{theorem}{Theorem}[section]
\newtheorem{definition}{Definition}[subsection]
\usepackage{listings}
\usepackage{caption}
\newtheorem{conj}{Conjecture}

\newtheorem{notation}{Notation}
\newtheorem{corollary}{Corollary}[theorem]
\newtheorem{def-lemma}{Definition-Lemma}
\newtheorem{rem}[theorem]{Remark}
\newtheorem*{cor*}{Corollary}
\newtheorem{lemma}[theorem]{Lemma}

\newtheorem*{rem*}{Remark}
\usepackage{verbatim}
\newtheorem*{lemma*}{Lemma}

\newtheorem*{claim*}{Claim}

\tikzset{
	symbol/.style={
		draw=none,
		every to/.append style={
			edge node={node [sloped, allow upside down, auto=false]{$#1$}}}
	}
}

\theoremstyle{definition}

\newcommand{\hilb}{\mathrm{hilb}_{\mathrm{ser}}}
\newcommand{\p}{\mathbb{P}}
\newcommand{\ri}{\rightarrow}
\newcommand{\J}{\mathscr{J}}
\newcommand{\jac}{\mathcal{J}}

\newcommand{\s}{\sigma}
\newcommand{\x}{\mathcal{X}}
\newcommand{\y}{\mathcal{Y}}
\newcommand{\m}{\mathcal{M}}

\newcommand{\h}{\mathcal{H}}

\newcommand{\U}{\mathscr{U}}
\newcommand{\ve}{\mathcal{V}ec}

\newcommand{\mcalu}{\mathcal{U}}

\newcommand{\E}{\mathscr{E}}
\newcommand{\e}{\mathcal{E}}
\newcommand{\Q}{\mathbb{Q}}
\newcommand{\C}{\mathscr{C}}
\newcommand{\mC}{\mathcal{C}}
\newcommand{\gm}{\mathbb{G}_m}
\newcommand{\ka}{\kappa}

\newcommand{\T}{\Theta}

\newcommand{\lb}{\mathcal{L}}

\newcommand{\spec}{\operatorname{Spec}}
\newcommand{\ext}{Filt}
\newcommand{\Hom}{\operatorname{Hom}}

\newcommand{\hr}{\hookrightarrow}

\usepackage{stackengine}
\newcommand{\hollowslash}{\setbox0=\hbox{/}\def\holwd{3pt}%
  \stackengine{-.3pt}{/}{\rlap{\kern.5pt\rule{\holwd}{0.4pt}}}{O}{r}{F}{F}{S}%
  \kern\dimexpr\holwd-\wd0-.2pt\relax%
  \stackengine{-.4pt}{/}{\llap{\rule{\holwd}{0.4pt}\kern-0.1pt}}{U}{l}{F}{F}{S}%
}

\title[]{\MakeUppercase{Cycles on the universal moduli of rank 2 bundles over the hyperelliptic locus} }
\author{Soumik Ghosh} 
\address{Department of Mathematics, Yale University, New Haven, CT 06511, USA} 
\email{soumik.ghosh@yale.edu}
\author{Shubham Saha}
\address{Department of Mathematics, University of California, San Diego, La Jolla, CA 92093, USA}
\email{shsaha@ucsd.edu}

\usepackage{microtype}
\begin{document}
	\begin{abstract}
We establish a strong Franchetta property for the universal moduli stack of rank two vector bundles over smooth hyperelliptic curves and identify the image of the Chow ring in the cohomology of a fiber via restriction and the cycle class map as precisely the monodromy invariant cohomology classes. More precisely, for every genus $g\geq 2$, degree $d$, and hyperelliptic curve $C$ of genus $g$, we show that the restriction and the cycle class map identifies the Chow ring of the universal moduli stack with the subring of $H^* (Bun_{2,d}(C),\mathbb Q)$ generated by twisted $\kappa$-classes, which coincides with $H^{*}(Bun_{2,d}(C), \mathbb Q)^{\pi_1(\mathcal H_g,[C])}$. We use this viewpoint to give an alternative proof to a theorem of Larson in \cite{larson-24}, and verify a conjecture of the second author in the rank $2$ case. We also give a general strategy for the higher rank cases and get some partial results in higher ranks and low genus. Finally, in Appendix \ref{MCK Appendix}, we establish the existence of multiplicative Chow-Kunneth decompositions for the moduli space of rank $2$, odd degree semi-stable bundles over low genus hyper-elliptic curves using the results developed in this paper. 
	\end{abstract}
	\subjclass[2020]{14C15, 14C17, 14D20, 14H45}
	\maketitle
    \thispagestyle{headings}
	\vspace{-1em}
    \setcounter{tocdepth}{1} 
\tableofcontents
\setcounter{section}{0}
\section*{Introduction}\label{intro}
Understanding algebraic cycles on moduli spaces of vector bundles over curves is a classical problem in algebraic geometry. While the cohomology of these moduli spaces are comparatively well-understood, describing their Chow groups is substantially more subtle. One approach is to study the Chow groups of the universal moduli spaces, where the variation of the underlying curve imposes additional rational equivalences on algebraic cycles. This leads naturally to Franchetta-type questions: to what extent are cycles on a universal family determined by their cohomology classes on individual fibers. The following conjecture was made by the second author in \cite{shubham-r2-ss}.  
\begin{conj}[\cite{shubham-r2-ss}]\label{conj:hyp}
	Let $\ve^\h({r,d,g})$ be the universal moduli stack of vector bundles of rank $r$, degree $d$ over hyperelliptic curves of genus $g$.
For any $g\geq 2, d\in\Z$ and a hyperelliptic curve $C$ of genus $g$, the composition of restriction and cycle-class map 
	$$A^*(\ve^\h({r,d,g}))\ri H^*(Bun_{r,d}(C))$$ is an isomorphism onto $H^*_{\ka}$ where $H^*_{\ka}$ is the cohomology subring generated by twisted-kappa classes.
\end{conj}
\noindent In loc. cit. the second author proves Conjecture A for rank $2$ bundles on genus $2$ curves with arbitrary degree. In this paper, we study this question for rank $2$ bundles over smooth hyperelliptic curves of any genus $g\geq 2$. 
Our main result is the following. \begin{theorem} \label{main theorem}
Conjecture \ref{conj:hyp} is true for $r=2$, any $d\in \mathbb Z$ and $g\geq 2$.     
\end{theorem}
\noindent We note that Conjecture \ref{conj:hyp} is false over $\mathcal M_g$ for $g\geq 3$ because of the Gross-Schoen cycle (see Remark \ref{Gross-Schoen}). We also obtain an explicit description of the rational Chow ring of $\mathcal Vec^{\mathcal H}(2,d,g)$ when $d$ is odd in terms of generators and relations among twisted $\kappa$ classes. In particular, using a monodromy argument, we obtain the following theorem.
\begin{theorem} \label{global Chow = monodromy invariant of a fiber}
    Let $\pi: \mathcal Vec^\mathcal H(2,d,g) \rightarrow \mathcal H_g$ be the universal stack. Fix a hyper-elliptic curve $C \in \mathcal H_g$. Then we have the isomorphism of graded rings \begin{equation*} 
        A^* (\mathcal Vec^\mathcal H(2,d,g)) \xrightarrow{\simeq } H^{2*} (Bun_{2,d}(C),\mathbb Q)^{\pi_1(\mathcal H_g, [C])}.
    \end{equation*}
\end{theorem}
\noindent The central result of \cite{larson-24} may be re-stated as follows.
\begin{theorem*}[\cite{larson-24}]
    The composition of restriction and the cycle class map gives the following isomorphism of graded rings for any hyper-elliptic curve $C$ of genus $g$,
    \begin{equation*}
        A^* (\mathcal Vec^\mathcal H(1,d,g)) \xrightarrow{\simeq} H^{2*} (Bun_{1,d}(C),\mathbb Q)^{\pi_1(\mathcal H_g, [C])} .
    \end{equation*}
\end{theorem*}\noindent So Theorem \ref{global Chow = monodromy invariant of a fiber} is a natural extension of Larson's result and we give a relatively simpler proof of Larson's result (see Theorem \ref{simpler proof of Larson}).  

A central ingredient in our argument is a relative version of the motivic decomposition for moduli spaces of stable bundles as in \cite{hoskins}. We extend this decomposition to families of curves over a smooth base. More precisely, let $\pi: \mathcal C \rightarrow B$ be a smooth projective family of curves of genus $g\geq 2$, where $B$ is smooth. $U_{\mathcal C/B}(2,d)$ is the relative moduli space of stable rank $2$ bundles of odd degree $d$. Write $J(\mathcal C/B)$ for the relative Jacobian and $\mathcal C_B^{(n)}$ for the relative symmetric powers. We prove the following result in Section 2 (see Theorem \ref{thm:decomp-rel}).

\begin{theorem} \label{motivic decomposition over a base}
      Let $U_B := U_{ \mC/B} (2,d, g)$ be the relative moduli space of rank two stable bundles of degree $d$ with $d$ odd. 
	We have the following isomorphism in $\mathrm M(B)$, the category of relative smooth Chow motives: 
	\begin{equation*}\label{relativeisomorphism}
        \mathfrak h_B (U_B)\simeq \mathfrak h_B (Jac (\mathcal C/B)) \otimes  \left (\mathfrak h _B([\mC_B^{(g-1)}] (g-1) \oplus  \bigoplus_{i=0}^{g-2} \mathfrak h_B([\mC_B^{(i)}/S_i]) \otimes (\mathbb Q(i) \oplus \mathbb Q(3g-3-2i)  )\right ).
    \end{equation*} 
\end{theorem}
\noindent We apply Theorem \ref{motivic decomposition over a base} to the universal curve over the hyperelliptic locus and results from \cite{mehdi} to establish a strong Franchetta Property (see Definition \ref{SFP-defn}) for $d$ odd. The identification of the image of the cycle class map on the cohomology of the fibers is taken up in Section \ref{sec:conj-hyp}. We give two proofs of surjectivity. In the first proof, we use representation theory of the symplectic group $\operatorname{Sp}_{2g}$ to compute dimensions of $A^i (\jac^d_{2,g}\times_{\h_g}\jac^{d'}_{2,g})$ and $A^i(\mC_{2,g}^{(n)}\times_{\h_g}\jac^g_{2,g})_{1\leq n\leq g+6}$ and use a monodromy argument to identify the image. This also enables us to get an explicit Hilbert series for $\mcalu^\h(2,d,g))$ when $d$ is odd, giving another numerical proof of surjectivity using the strong Franchetta property. 
    The case $d$ is even is taken up in Section \ref{sec:even-deg}. We study the moduli stack $\mathcal H^w_{g,n}$ of hyperelliptic curves with $n$ distinct marked points and an additional Weierstrass point which can possibly coincide with one of the $n$ -markings. We show the following, after a modification of the arguments in \cite{cl-hyp} (see Theorem \ref{thm:weier-pull}).
    \begin{theorem} \label{additional Weirestrass point stack}
          For all $g\geq 2, 0\leq n\leq 2g+4$, the pullback map $A^*(\h_{g,n})\ri A^*(\h_{g,n}^w)$ is an isomorphism.
    \end{theorem}
\noindent This enables us to prove the strong Franchetta property and Conjecture \ref{conj:hyp} for $r=2$ and all even $ d$. In Section 8, we give a strategy for higher ranks. In particular, we prove the following. 

\begin{theorem} \label{higher rank}
    Let $\mathcal U^\mathcal H(r,d,g)$ denote the moduli stack of rank $r$ and degree $d$ semi-stable bundles over $\mathcal H_g$ with $\gcd(r,d)=1$. Then the family $\mathcal U^\mathcal H(r,d,g) \rightarrow \mathcal H_g$ satisfies the strong Franchetta property whenever $rg\leq 2g+6$. 
\end{theorem}
\noindent Finally, in Appendix \ref{MCK Appendix}, we study the existence of multiplicative Chow-Kunneth (MCK) decompositions for the moduli space of rank $2$ bundles over a hyper-elliptic curve. The following is our main result about the existence of MCK decompositions. 
\begin{theorem}\label{MCK Theorem}
    If $C$ is a hyper-elliptic curve of genus $g\leq 7$, then $U_C(2,\mathcal L)$ and $U_C(2, d)$ admit multiplicative Chow-Kunneth decompositions where $\mathcal L$ is a line bundle of odd degree and $d$ is an odd integer.  
\end{theorem}
\subsection*{Acknowledgements} 
SG would like to thank his PhD advisor, Junliang Shen, for his interest in the project, and David Bai and David Fang for helpful discussions. SS would like to thank his PhD advisor, Elham Izadi, for her guidance and support, and Qizheng Yin and Samir Canning for helpful discussions.

\subsection*{Acknowledgement of AI use}
The authors acknowledge the use of AI in the preparation of this manuscript. Google Gemini 3.1 Deep Think model via Perplexity was used to assist with the proofs of Theorems \ref{thm:ext-hilb}, \ref{thm:hilb-stable-denom}, and the code in Figure \ref{hilb:code}. The strategy of Theorem \ref{thm:ext-hilb} given by AI was flawed. However, the use of Equation \ref{AI assist} in the flawed 'proof' helped the authors to simplify the proof of Theorem \ref{thm:jxj-rel}. The proof of Theorem \ref{thm:ext-hilb} in this pre-print is the authors' own. For Theorem \ref{thm:hilb-stable-denom}, the idea of using residues and contour integration to get a closed form for the Hilbert series came from AI. The authors subsequently did the computations themselves. The authors further give a proof of the main theorems independent of Theorem \ref{thm:hilb-stable-denom}. The authors also acknowledge the usage of Open AI's ChatGPT Plus to do an extensive literature search on the representation theory of the symplectic group used in the proofs of Section \ref{sec:conj-hyp}. All the mathematical arguments and proofs appearing in the paper were written and verified by the authors. The authors take sole responsibility for all statements and arguments.

\section*{Notation}\label{notation}
	\noindent All schemes in this paper are defined over the field of complex numbers. 
	All Chow rings and groups are taken with rational coefficients. We use the subspace convention for projective bundles. 
	
	We denote the locus of hyperelliptic curves of genus $g$ with $n$-marked points with an additional marked Weierstrass point as $\mathcal{H}^w_{g,n}\hr \h_{g,1}\times_{\h_g}\h_{g,n}$, and put $\mC^w_{2,g,n}$ as the universal curve over $\mathcal{H}^w_{g,n}$ with the Weierstrass section $\{\s_{g,i,n}: \mathcal{H}^w_{g,n}\ri \mC^w_{2,g,n}\}$.  

    Analogous to \cite{cl-hyp}, we set $D_{i,j}:= \left \{ (C,p_1,\dots, p_n)\in \mathcal H_{g,n}\mid p_i = \bar p_j \right \}$ which is the locus where $p_i$ and $p_j$ are conjugate under the hyper-elliptic involution. 
	
	For a fibre product $\tensor[_{i\in I}]{\prod}{_S}X_i$ over a scheme $S$, and a subset $J\subset I$, we denote the projection onto the sub-product $\tensor[_{j\in J}]{\prod}{_S}X_j$ by $p_J$. 

    By Hilbert Series of a graded $\Q$-vector space $V^\bullet$, we refer to the series $\hilb(V^\bullet) := \sum_{i\geq 0}\dim_\Q(V^i)T^i$ where $T$ is an indeterminate. For a stack $\chi$, we refer to the series $\hilb(\chi):= \sum_{i\geq 0}T^i\dim A_{\dim \chi - i}(\chi) $.

For a stack $\mathcal{X}$, $H^*_a(\mathcal{X})$ denotes the algebraic cohomology subring of $H^*(\mathcal{X})$, i.e., the image of the cycle class map $A^*(\mathcal{X})\ri H^*(\mathcal{X})$.

We denote the universal Jacobian over hyperelliptic curves of genus $g$ with degree $d$ as $\J^d_{2,g}$ and its rigidification by $\jac^d_{2,g}$. The associated $\mathbb{G}_m$-rigidification map is denoted by $\nu_{1,d}:\J^d_{2,g}\ri\jac^d_{2,g}$.
\section{Preliminaries}

\subsection{$\mathbb{G}_m$-gerbes}
For technical simplicity, we shall use the following notion of a torsion-$\gm$-gerbe. 
\begin{definition}[\cite{min-shen}]
    A \textit{torsion-$\mathbb{G}_m$-gerbe} is a gerbe $\mathcal{Y}\ri \mathcal{X}$, banded by $\mathbb{G}_m$, such that the associated cohomology (Brauer-Grothendieck) class is torsion. That is, $[\mathcal{Y}/\mathcal{X}]\in H^2(\mathcal{X}_{\text{\'et}},\mathbb{G}_m)_{tors}$.
\end{definition}
Note that any $\gm$-gerbe defined over a regular, integral, and noetherian scheme is a torsion-$\gm$-gerbe \cite[Corollaire 1.8]{grothendieck}. The rigidification map $\J^d_{g}\ri \jac^d_{g}$ is also a torsion-$\gm$-gerbe by \cite[Theorem 6.4]{melo}.
\begin{theorem}[\cite{larson-24}]\label{thm:gerbe-ring}
    Let $\mathcal{Y}\ri \mathcal{X}$ be a torsion-$\gm$-gerbe. There exists $c_1\in A^1(\mathcal{Y})$ such that
    $$A^*(\mathcal{Y}) \simeq A^*(\mathcal{X})[c_1].$$
\end{theorem}
\begin{proof}
    Let $n\in \mathbb{N}$ be such that $n[\mathcal{Y}/\mathcal{X}] = 0$. Applying \cite[Lemma 3.1.1.8]{lieblich}, there exists a $\mu_n$-gerbe $\mathcal{Y}\ri\mathcal{Y}_n$ such that
$$\begin{tikzcd}
    \mathcal{Y}\ar[r]\ar[rd]&\mathcal{Y}_n\simeq \mathcal{X}\times B\gm\ar[d]\\
    &\mathcal{X}.
\end{tikzcd}$$
Therefore, we have $A^*(\mathcal{Y})\simeq A^*(\mathcal{Y}_n)$. Applying the classical equivariant-Chow argument from \cite{edidin-graham}, we have 
$$A^*(\mathcal{X}\times B\gm)\simeq A^*(\mathcal{X})\otimes A^*(B\gm)\simeq A^*(\mathcal{X})[c_1].$$
\end{proof}
\subsection{Twisted Sheaves}\label{sec:twisted}
We follow Lieblich's notion of twisted sheaves \cite{lieblich-main}.  A \textit{twisted-sheaf} on a smooth DM stack $X$ is a coherent sheaf on $\mathcal{X}$ where 
$\mathcal{X}\ri X$ is a gerbe banded by a finite abelian group. To emphasize this distinction between twisted sheaves and coherent sheaves, we shall use  coherent sheaves and \textit{honest coherent sheaves} interchangeably. 

For a torsion Brauer class $\alpha\in Br(X)$ such that $n\alpha = 0$, let $\mathscr{X}$ be the associated $\gm$-gerbe. We can descend honest coherent sheaves on $\mathscr{X}$ using the following theorem.
\begin{theorem}[\cite{lieblich-main}]\label{thm:twisted}
    Suppose we have a coherent sheaf $\mathcal{F}\ri \mathscr{X}$, then $\mathcal{F}$ descends to $X$ as a twisted-sheaf.
\end{theorem}
\begin{proof}
We use the Kummer sequence and the long exact sequence in \'etale cohomology
$$1\ri \mu_n\ri \gm\xrightarrow{(-)^n}\gm\ri 1.$$
We have the exact sequence $$H^1(X,\gm)\ri H^2(X,\mu_n)\xrightarrow{i} H^2(X,\gm)\ri H^2(X,\gm).$$
Since $n\alpha=0$, we have $\alpha\in Im(i)$. That is, there exists $ \beta\in H^2(X,\mu_n)$ such that $i(\beta)=\alpha$. The class $\beta$ corresponds to a $\mu_n$-gerbe $\mathcal{X}\ri X$. 

The restriction of $\alpha$ to $\mathcal{X}$ is trivial \cite[Theorem 1.2]{min-shen} since $0= \beta|_{\mathcal{X}}$ maps to $\alpha|_{\mathcal{X}}$ under $i|_{\mathcal{X}}$. Hence, we have the Cartesian diagram
$$\begin{tikzcd}
    \mathscr{X}\times B\gm \ar[r,"\alpha|_{\mathcal{X}}"]\ar[d]\ar[dr, phantom, "\square"]&\mathscr{X}\ar[d]\\
    \mathcal{X}\ar[r]\ar[u, "s", bend left = 60]&X.
\end{tikzcd}$$
The pullback $(\alpha|_{\mathcal{X}}\circ s)^*\mathcal{F}$ is an honest coherent sheaf on $\mathcal{X}$.
Hence, we have a twisted sheaf on $X$ with the claimed property.
\end{proof}
\begin{definition}
    We shall refer to these descended twisted-sheaves on $X$ as $\alpha$-twisted sheaves.
\end{definition}

\subsection{The Strong Franchetta Property}
For a morphism $f: \mathcal{F}\ri \mathcal{Z}$ between smooth, finite-type algebraic stacks, we say that $f$ satisfies the Franchetta property if:\\
For every $\Gamma\in A^*(\mathcal{F})$ such that $\Gamma|_{\mathcal{F}_z}$ is trivial in $H^*(\mathcal{F}_z)$ for a very general geometric point $z\in |\mathcal{Z}|$, the restriction $\Gamma|_{\mathcal{F}_z}$ is trivial in $A^*(\mathcal{F}_z)$ for all $z\in |\mathcal{Z}|$.

This property has been studied in \cite{laterveer-franchetta, lifranchetta}. We refer to these sources for more details.
For our purposes, we define a variant of the property below.
\begin{definition}\label{SFP-defn}
  We say that a family $\phi: \mathcal{X}\ri \mathcal{S}$ satisfies the \textit{strong Franchetta property} (SFP) if the map $A^*(\mathcal{X})\ri H^*({X}_s)$ is  injective for all fibers $X_s$ over geometric points $s:\spec \mathbb{C}\ri \mathcal{S}$.
\end{definition}

\begin{theorem}\label{thm:gm-sfp}
   Let $\mathcal{X}\ri S$ be a family satisfying SFP, and let $\mathcal{Y}\ri\mathcal{X}$ be a torsion-$\mathbb{G}_m$-gerbe. Then $\mathcal{Y}\ri S$ satisfies SFP. 
\end{theorem}
\begin{proof}
    For a geometric point $s\in |S|$, and the fibers $Y_s,,X_s$ of $\mathcal{Y},\mathcal{X}$ over it respectively, we have  
    $$Y_s\ri X_s$$ is a $\gm$-gerbe. Applying Theorem \ref{thm:gerbe-ring} and \cite[Lemma 5.13]{MR3013029}, we have that there exists $c_1\in A^1(\mathcal{Y})$ such that $$A^*(\mathcal{Y})\simeq A^*(\mathcal{X})[c_1], H^*(Y_s)\simeq H^*(X_s)[c_1|_{Y_s}].$$
    SFP of $\mathcal{X}\ri S$ gives the following diagram:
$$\begin{tikzcd}
A^*(\mathcal{X})[c_1]\ar[r, dashed]& H^*(X_s)[c_1|_{Y_s}]\\
A^*(\mathcal{X})\ar[u, hook]\ar[r, hook]& H^*(X_s)\ar[u,hook]
\end{tikzcd}
$$
Injectivity of the dashed horizontal arrow is immediate from the diagram. Consequently, we have SFP for $\mathcal{Y}\ri S$. 
\end{proof}

\begin{theorem}\label{thm:inj-sfp}
    Let $\mathcal{X}\xrightarrow{f}\mathcal{Y}$ be a morphism of families over $S$, with $\mathcal{X},\mathcal{Y}$ smooth such that $f^*:A^*(\mathcal{Y})\ri A^*(\mathcal{X})$ is injective. The family $\mathcal{Y}\ri S$ satisfies SFP if we have the same for $\mathcal{X}\ri S$.
\end{theorem}
\begin{proof}
For fibers $Y_s,X_s$ of $\mathcal{Y},\mathcal{X}$ respectively, over a geometric point $s\in |S|$, we have the following diagram due to SFP of $\mathcal{X}\ri S$:
$$\begin{tikzcd}
A^*(\mathcal{Y})\ar[r, dashed]\ar[d,"f^*", hook]& H^*(Y_s)\ar[d]\\
A^*(\mathcal{X})\ar[r, hook]& H^*(X_s).
\end{tikzcd}
$$
Injectivity of the dashed arrow follows from the diagram, which shows SFP for $\mathcal{Y}\ri S$.
\end{proof}

The following theorem follows from a straightforward application of the projective bundle formula in Chow and cohomology.
\begin{theorem}\label{thm:sfp-proj}
Let $\mathcal{Y}\ri \mathcal{X}$ be a projective bundle and $\mathcal{X}\ri S$ be a family. 
We have that $\mathcal{X}\ri S$ satisfies SFP if and only if $\mathcal{Y}\ri S$ satisfies SFP.
\end{theorem}

\section{A motivic decomposition over families}\label{sec:motive-families}

\subsection{Chow Motives over a Base}

Let $B$ be a non-singular variety over $\mathbb C$. Following \cite{MR3851129}, we denote the category of relative smooth Chow motives over $B$ with $\mathbb Q$-coefficients as $\mathrm{M}(B)$. 
The objects in $\mathrm{M}(B)$ are $(X\xrightarrow{f}B, \mathfrak p, m)$ where $X\xrightarrow{f}B$ is smooth projective, $\mathfrak p \in \operatorname{Corr}_B ^0(X,X)$ is a projector and $m\in \mathbb Z$ records the Tate twists. Further, we have \begin{equation}
\operatorname{Hom}_{\mathrm{M}(B)} \left ( (X, \mathfrak p , m), (Y, \mathfrak q, n) \right ) = \mathfrak q\circ \operatorname{Corr}_B^{n-m}\circ \mathfrak p.
\end{equation}
Further, we denote by $\mathrm A(B)$ the full pseudo-abelian tensor category generated additively and tensorially by motives of relative curves in $\mathrm{M}(B)$. 

\begin{theorem}\cite[Corollary C]{MR3851129}\label{relative-nilpotence-theorem}
    Let $\mathfrak h $ be a motive in $\mathrm{A}(B)$, let $\eta$ be the generic point of $B$ with residue field $\mathbb C(\eta)$, let $L/{\mathbb C(\eta)}$ be a field extension, and set $\xi:=\operatorname{Spec}L$. Any endomorphism of a motive in $\mathrm A(B)$ whose pull-back to $\xi$ is nilpotent modulo numerical equivalence is nilpotent modulo rational equivalence over the base scheme $B$. 
\end{theorem}
\begin{corollary}\label{generic-isomorphism}
    Suppose we have two motives $\mathfrak h_1, \mathfrak h_2$ in $\mathrm A(B)$ and $\Gamma\in \operatorname{Hom}_{\mathrm{M}(B)}(\mathfrak h_1, \mathfrak h_2), \Psi\in \operatorname{Hom}_{\mathrm{M}(B)}(\mathfrak h_2, \mathfrak h_1)$
	 be two correspondences such that for all $b\in B$ closed points, \\
	 $\Gamma_b \in \operatorname{Hom}_{\mathrm{M}(\mathbb C)}\left((\mathfrak h_1)_b, (\mathfrak h_2)_b\right)$, $\Psi_b\in \operatorname{Hom}_{\mathrm{M}(\mathbb C)} \left ( (\mathfrak h_2)_b, (\mathfrak h_1)_b \right )$ satisfies $\Gamma_b \circ \Psi_b =Id_{(\mathfrak h_2)_b}$ and $\Psi_b\circ \Gamma_b=Id_{(\mathfrak h_1)_b}$. Then $\Gamma\in \operatorname{Hom}_{\mathrm{M}(B)}(\mathfrak h_1, \mathfrak h_2)$ is an isomorphism. 
\end{corollary}
\begin{proof}
    Consider the correspondences $E_1= Id_{\mathfrak h_1}-\Psi\circ \Gamma$ and $E_2= Id_{\mathfrak h_2}-\Gamma\circ \Psi$. 
	Then we have by condition that $E_1|_b=0, E_2|_b =0$ for all $b \in B$ closed points. 
	Applying \cite[Lemma 2.1.2]{voisin-cycles}, it follows that $(E_1)_\eta =0, (E_2)_\eta=0$ where $\eta$ is the generic point of $B$ and hence by Theorem \ref{relative-nilpotence-theorem}, we see that $E_1\in \operatorname{End}_{\mathrm{M}(B)}(\mathfrak h_1), E_2\in \operatorname{End}_{\mathrm{M}(B)}(\mathfrak h_2)$ are nilpotent endomorphisms. Therefore $\Psi\circ \Gamma$
 and $\Gamma\circ \Psi$ are isomorphisms. It follows that $\Gamma$ is an isomorphism. \end{proof}
 
\subsection{Motive of Moduli of Stable Bundles over a Base}\label{Motives of Bundles on Curves over a Base}
  Let $\mathcal C \rightarrow B$ be such that for any closed point $b\in B$, the fiber $C_b$ is a smooth, proper curve of genus $g\geq 2$. Let $ U_{\mathcal C/B}(r,d,g)$ be the relative coarse moduli space of stable bundles of rank $r$ and degree $d$ with $r,d$ coprime. 
  
  We have $ U_{\mathcal{C}/B} (r,d,g)\ri B$ is a smooth proper morphism of of  relative dimension $m:=r^2(g-1)+1$. For notational convenience, set $U:=U_{\mathcal C/B}(r,d,g)$. We shall use a relative version of Beauville's diagonal argument as in \cite{MR1351502} to show $\mathfrak h_B(U)$ is Kimura finite-dimensional. We refer to \cite{MR2107443} and \cite{MR3851129} for the necessary definitions. For our case, we shall only need that the relative motive is in $\mathrm A(B)$. 

\begin{lemma}\label{lem2.1}
    Let $\Delta_{U/B}$ be the image of the embedding $U\hr U\times_B U$. We have \begin{equation*}[\Delta_{U/B}] \in \Q. c_{m} \left (-\mathscr Ext^!_{\pi_{13}} \left (\pi_{12}^* \mathcal E', \pi_{23}^*\mathcal E''  \right ) \right )\end{equation*} in $A^m \left (U\times_B U \right )$ where the morphisms $\pi_{ij}$'s are given by the projections in the diagram  $$\begin{tikzcd}&  U \times_B \mathcal C \times_B  U \arrow{d}{\pi_{13}}\arrow{rd}{\pi_{23}}\arrow{ld}{\pi_{12}} \\  U \times_B \mathcal C &  U\times_B  U & \mathcal C\times_B  U\end{tikzcd}. $$
\end{lemma}
\begin{proof}
Let $\alpha$ be the Brauer class of the $\gm$-gerbe $\U\ri U$ where $\U$ is the unrigidified moduli stack of fiber-wise stable, rank $r$ bundles with relative degree $d$ for the family $\mC\ri B$. We note that $\alpha \in \operatorname{Br}(U)$ is torsion by \cite[Theorem B]{brauer}.
Let $\mathcal E', \mathcal{E''}$ be twisted universal sheaves on $U\times_B \mathcal C, \mC\times_B U$ respectively, defined using $\alpha$ and the universal sheaves on $\U\times_B\mC, \mC\times_B\U $ respectively by Theorem \ref{thm:twisted}. 

Let $ N\in \mathbb N$ be the order of $\alpha$ and $\mathcal G\xrightarrow{g} U$ be the associated $\mu_N$-gerbe.
   It follows that $\mathcal E'$ is $q_1^*\alpha$-twisted and $\mathcal E''$ is $r_2^* \alpha$-twisted where $q_1: U\times_B \mathcal C \rightarrow U$, $r_2: \mathcal C\times_B U\rightarrow U$ are the projections. 

   We denote the associated honest coherent sheaves on $\mathcal{G}\times_B\mC, \ \mC\times_B \mathcal{G}$ as $\tilde{\e'}, \ \tilde{\e''}$ respectively with projections
   $\tilde{\pi}_{12}:\mathcal G\times_B \mathcal C \times_B \mathcal G\ri \mathcal G\times_B \mathcal C , \ \tilde{\pi}_{23}: \mathcal G\times_B \mathcal C \times_B \mathcal G\ri  \mathcal C \times_B \mathcal G$ onto (first, second) and (second, third) factors respectively.
   
   
 The pull-backs $\tilde{\pi}_{12}^*\tilde{\e'}, \tilde{\pi}_{23}^*\tilde{\e''}$, are honest coherent sheaves on $\mathcal G\times_B \mathcal C \times_B \mathcal G$.   Consequently,  $R {{\tilde \pi}_{13*}}R\mathscr Hom \left (R{\tilde \pi}_{12}^* \tilde {\mathcal E'} , R\tilde \pi_{23}^*\tilde{\mathcal E''}  \right )$ can be represented on $\mathcal G\times_B \mathcal G$ by a two-term complex $K^0\rightarrow K^1$. 
   
   For a geometric point $(\xi, \eta)\in \mathcal G\times_B \mathcal G$ mapping to $(x, y)\in U\times_B U$, the fiber of the complex computes $RHom (E_x, E_y)$ and $R^0Hom (E_x, E_y) =\begin{cases}
       \mathbb C \text{ if }x=y \\ 0 \text{ if }x\neq y
   \end{cases}$. Using the Porteous formula for smooth DM stacks , we get
   $$[\mathcal G \times_U \mathcal G] \in \Q.c_m\left (-R {{\tilde \pi}_{13*}}R\mathscr Hom \left (R{\tilde \pi}_{12}^* \tilde {\mathcal E'} , R\tilde \pi_{23}^*\tilde{\mathcal E''}  \right )\right )\in A^*(\mathcal{G}\times_B\mathcal{G}).$$
  Since $(g\times_B g)_* : A^*(\mathcal G\times_B \mathcal G) \rightarrow A^*(U\times_B U)$ is an isomorphism, we have \begin{equation*}
       [\Delta_{U/B}] \in \Q. c_{m} \left (-\mathscr Ext^!_{\pi_{13}} \left (\pi_{12}^* \mathcal E', \pi_{23}^*\mathcal E''  \right ) \right ).
   \end{equation*}
 \end{proof}
  \begin{lemma}\label{lem2.2}
    Let $M,S,S',M'$ be smooth proper $B$-varieties. If $\alpha \in \operatorname{Corr}_B(M,S), \beta \in \operatorname{Corr}_B(S, M'),\\ \alpha'\in \operatorname{Corr}_B(M,S'), \beta'\in \operatorname{Corr}_B(S',M')$, then 
	\begin{equation*}
        (\beta'\circ \alpha')\cdot (\beta\circ \alpha)= \gamma \circ \delta 
    \end{equation*}
	 for some $\gamma \in \operatorname{Corr}_B(S\times_B S', M'), \delta \in \operatorname{Corr}_B(M, S\times_B S')$.
\end{lemma}
\begin{proof}
    It follows from a straightforward computation that \begin{equation*}
        \begin{aligned}
         (\beta'\circ \alpha')\cdot (\beta\circ \alpha) =  \left ({P'_{13}}^* \beta\cdot {P'_{23}}^* \beta'  \right ) \circ \left ({P_{12}}^* \alpha\cdot {P_{13}}^* \alpha' \right ) \text{ where }
        \end{aligned}
    \end{equation*} 
    $$\begin{tikzcd}
       & S\times_B S' \times_B M'\arrow{rd}[swap]{P'_{23}}\arrow{ld}{P'_{13}} & \\ S\times_B M' & &S'\times_B M'
    \end{tikzcd} ,\begin{tikzcd}
       & M\times_B S\times_B S' \arrow{rd}[swap]{P_{13}}\arrow{ld}{P_{12}} & \\ M\times_B S & &M\times_B S'
    \end{tikzcd}$$ are the projection morphisms.  
\end{proof}
 By Grothendieck-Riemann-Roch, we have \begin{equation}\label{eq2}
    ch \left ( \mathscr Ext^!_{\pi_{13}} \left (\pi_{12}^* \mathcal E', \pi_{23}^* \mathcal E'' \right ) \right ) = \left ( ch (\mathcal E''). q_1^* \sqrt{td_\mC}\right )\circ \left ( ch (\mathcal E')^\vee. p_2^* \sqrt{td_\mC}\right )
\end{equation} 
where $q_1: \mC\times_B  U\rightarrow \mC, p_2:  U\times_B \mC \rightarrow  \mC$ are the projections, and $\sqrt{td_\mC}$ is the formal square root of the unit $td_\mC\in A^*(\mC)$. 
For notational convenience, we shall denote $\times_B^n\mC$ by $\mC^n$ for all $n\in\mathbb{N}$.

\begin{lemma}\label{lem:ck-gamma-delta}
 For any $k\geq 1$, there exist $\delta_{\alpha}\in \operatorname{Corr}_B( U, \mC^{k_\alpha}), \gamma_{\alpha }\in \operatorname{Corr}_B(\mC^{k_\alpha},  U)$ with $k_\alpha \leq k$ such that
   \begin{equation*}
        c_k \left (-\mathscr Ext^!_{\pi_{13}} \left (\pi_{12}^* \mathcal E', \pi_{23}^*\mathcal E''  \right ) \right ) =\sum_{\alpha }\gamma_\alpha \circ \delta_\alpha.
    \end{equation*} 
\end{lemma}
\begin{proof}
    By (\ref{eq2}), we have 
	$$ ch \left ( \mathscr Ext^!_{\pi_{13}} \left (\pi_{12}^* \mathcal E', \pi_{23}^* \mathcal E'' \right ) \right ) = \gamma \circ \delta $$ 
	for $\gamma \in \operatorname{Corr}_B(\mC,  U)$, $\delta\in \operatorname{Corr}( U, \mC)$. 
	We have
	\begin{equation*}
        ch = r+ c_1 + \left(\frac{1}{2}c_1^2 -c_2\right)+\cdots = \sum_{j=0}^\infty P_j 
    \end{equation*} 
	where $P_j$ is a homogeneous polynomial with $\mathbb Q$ coefficients in the Chern classes $c_1, c_2, \dots , c_j$ with $\deg c_k=k$ and the coefficient of $c_j$ in $P_j$ is non-zero.  
	A straightforward induction and Lemma \ref{lem2.2} show that  \begin{equation*}
        c_k \left ( -\mathscr Ext^!_{\pi_{13}} \left (\pi_{12}^* \mathcal E', \pi_{23}^*\mathcal E''  \right )\right ) =\sum_{\alpha }\gamma_\alpha \circ \delta_\alpha 
    \end{equation*} 
	where $\delta_{\alpha}\in \operatorname{Corr}_B( U, \mC^{k_\alpha}), \gamma_{\alpha }\in \operatorname{Corr}_B(\mC^{k_\alpha},  U)$ with $k_\alpha \leq k$. 
  
\end{proof} 

We present a relative version of \cite[Theorem 0.1]{MR4046978} below.
\begin{corollary}\label{U has abelian motive}
    We have $\mathfrak h_B (U)\in \mathrm A(B)$. 
\end{corollary}
\begin{proof}
Let $I\subset \operatorname{Corr}_B(U,U)$ denote the two-sided ideal generated by relative correspondences which factor through some relative power of $\mathcal C$. More precisely, we have \begin{equation*}
     I = \langle  \beta\circ \alpha \mid \alpha \in \operatorname{Corr}_B(U,\mC^k), \beta \in \operatorname{Corr}_B(\mC^k,U)\rangle .
 \end{equation*}
 We note that $I$ is closed under taking intersection products by Lemma \ref{lem2.2}. 
  Applying Lemma \ref{lem2.1} and Lemma \ref{lem:ck-gamma-delta}, there exist relative correspondences 
 $\delta_i \in \operatorname{Corr}_B(U,\mC^{k_i}), \gamma_i\in \operatorname{Corr}_B(\mC^{k_i},U) $, and $k_i\in \mathbb{N}$ such that 
 \begin{equation*}
     [\Delta_{U/B}] = \sum_i \gamma_i \circ \delta_i.
 \end{equation*}
 Viewing $\delta :=\oplus_i \delta_i$ as a morphism of motives \begin{equation*}
      \bigoplus_i \mathfrak h_B (\mathcal C^{k_i})(d_i-2k_i)\xrightarrow{\delta}\mathfrak h_B(U),
 \end{equation*} 
 we obtain the right inverse $\gamma:= \oplus_i \gamma_i $. Equivalently, we have a section of $\delta$
 \begin{equation*}
     \mathfrak h_B (U) \xrightarrow{\gamma} \bigoplus_i \mathfrak h_B (\mathcal C^{k_i})(d_i-2k_i).
 \end{equation*}
 Hence $\mathfrak h_B(U)$ is a direct summand of $\bigoplus_i \mathfrak h_B (\mathcal C^{k_i})(d_i-2k_i)$.
 Consequently, we have the claimed containment.
\end{proof}

\noindent Assume now that we have a genus $g$ curve fibration $\mathcal C \rightarrow B$ with a section $\sigma$. So there exists a line bundle $\mathcal L$ on $\mathcal C/B$ such that $L_b$ the restriction of $\mathcal L$ on $C_b$, a fiber over $b\in B$ gives a degree $1$ line bundle. We consider the family $U_B (r,\mathcal L, g)\rightarrow B$ which over every point parametrizes the semi-stable rank $r$ vector bundles with determinant $L_b$. A straightforward application of Beauville's diagonal trick as established in \cite[Proposition 4.1]{hoskins}, combined with the analogous arguments of Lemma \ref{lem2.1} and Lemma \ref{lem:ck-gamma-delta} gives us the following corollary.
\begin{corollary}
   We have  $\mathfrak h_B (U(r,\mathcal L, g) )\in \mathrm A (B)$. 
\end{corollary}

\subsection{A Relative Motivic Decomposition for Moduli of Bundles}\label{sec:rel-motive-decomp}
 Let $\mathcal C\rightarrow B$ be a family of smooth curves as above. Fix $d\in\Z$ odd for the rest of this section. The following is proved in \cite[Theorem 4.3]{hoskins}

 \begin{theorem}[\cite{hoskins}]\label{HoskinsMotivDecomp}
 For a smooth projective curve $C$ of genus $g$ over an arbitrary field $k$ of $\operatorname{char} 0$, let $U:= U_C(2,d)$ be the moduli space of stable bundles of rank $2$ and degree $d$. We have the following isomorphism in $M(k)$ 
     \begin{equation*}
         \mathfrak h(U) \simeq \mathfrak h(Jac (C)) \otimes \left (\mathfrak h (Sym^{g-1} C)(g-1) \oplus  \bigoplus_{i=0}^{g-2} \mathfrak h(Sym^i C) \otimes (\mathbb Q(i) \oplus \mathbb Q(3g-3-2i) ) \right )
     \end{equation*}
 \end{theorem}
\begin{proof} The proof in \cite{hoskins} goes through the existence of a line bundle of odd degree, which we want to avoid. So we indicate a sketch of the above decomposition following their argument. No novelty is claimed as to the proof.

\noindent It follows from \cite[Theorem 4.11]{del-bano} and the argument in \cite[Theorem 4.3]{hoskins} that the following relation holds in $K_0 (\operatorname{\mathrm M}_{num}(k))$: \begin{equation*}
    [\mathfrak h(U)] = [\mathfrak h(Jac(C))] \ \left ( [\mathfrak h(Sym^{g-1}C)]\mathbb L^{g-1} +\sum_{i=0}^{g-2} [\mathfrak h(Sym^i C) ] \left ( \mathbb L^i +\mathbb L^{3g-3-2i} \right )\right ). 
\end{equation*}
Moreover we know that $\mathfrak h(U) , \mathfrak h(Jac(C)) , \mathfrak h( C) ,\mathbb L \in \mathrm A(k)$ by \ref{U has abelian motive},  \cite[Theorem 4.5]{del-bano} (this is where we use $\operatorname{char} k=0$, as both \ref{U has abelian motive} and \cite[Theorem 4.5]{del-bano} require that hypothesis).  By the results in \cite{MR1150598}, we can lift this to an isomorphism of numerical motives and since both sides are of abelian type, we can further lift them to an isomorphism in $\mathrm M(k)$.    
\end{proof}
The objective of this section is to extend the above result for smooth curve fibrations $\mathcal C\rightarrow B$. 

\begin{theorem}\label{thm:decomp-rel}
    Let $U_B := U_{ \mC/B} (2,d, g)$ be the relative moduli space of rank two stable bundles of degree $d$. 
	We have the following isomorphism in $\mathrm M(B)$: 
	\begin{equation*}\label{relativeisomorphism}
        \mathfrak h_B (U_B)\simeq \mathfrak h_B (Jac (\mathcal C/B)) \otimes  \left (\mathfrak h _B([\times_B^{g-1} \mathcal C/S_{g-1}](g-1) \oplus  \bigoplus_{i=0}^{g-2} \mathfrak h_B([\times_B^i \mathcal C/S_i]) \otimes (\mathbb Q(i) \oplus \mathbb Q(3g-3-2i)  )\right )
    \end{equation*} 
\end{theorem}
\begin{proof}
    Let $\eta \in B$ be the generic point with residue field $k(\eta)$. Applying Theorem \ref{HoskinsMotivDecomp}, we have the following isomorphism in $\mathrm{M} (k(\eta))$:
	 \begin{equation*}
        \mathfrak h(U_\eta) \simeq \mathfrak h(Jac(\mathcal C_\eta)) \otimes \left (\mathfrak h (Sym^{g-1} \mathcal C_\eta)(g-1) \oplus  \bigoplus_{i=0}^{g-2} \mathfrak h(Sym^i \mathcal C_\eta) \otimes (\mathbb Q(i) \oplus \mathbb Q(3g-3-2i)  )\right ).    
	\end{equation*} 
        Let $\Gamma_\eta $ be a correspondence that induces an isomorphism above with inverse given by $\Psi_\eta$,
 let $\mathfrak M: = \mathfrak h_B (U)$, and let $\mathfrak N$ be the motive on the right of the Theorem statement. 
		
 Using the ``spread-out" argument by Voisin \cite{voisin-spread}, there exist
           $ \Gamma \in \Hom_{\mathrm M(B)} \left ( \mathfrak M, \mathfrak N\right ) , \Psi \in \operatorname{Hom}_{\mathrm M(B)} (\mathfrak N , \mathfrak M)$
		 such that $\Gamma$ specializes to $\Gamma_\eta$ and $\Psi$ specializes to $\Psi_\eta$ at $\eta$. 
		By specialization it follows that for every closed point $b\in B$, we have $\Gamma_b \circ \Psi_b=Id$ and $\Psi_b\circ \Gamma_b=Id$. Applying Corollary \ref{generic-isomorphism}, we have that $\Gamma$ is an isomorphism since $\Gamma_\eta, \Psi_\eta$ are inverses of each other and $\mathfrak M, \mathfrak N\in \mathrm{A}(B)$.   
\end{proof}

\section{Algebraic and tautological cohomology for bundles on a curve}

Let $\mathcal C \rightarrow B$ be as in \S\ref{Motives of Bundles on Curves over a Base}. We recall that we have the following equality from \ref{lem2.1} and \ref{lem:ck-gamma-delta} and \begin{equation}\label{diagonal as sum of composition}
    [\Delta_{U/B}] = \sum_i \gamma_i \circ \delta_i
\end{equation} 
for $\delta_i \in\operatorname{Corr}_B(U, \mathcal C^{k_i})$  and $\gamma_i \in \operatorname{Corr}_B (\mathcal C^{k_i} , U) $. 
Further, it follows from the proof of Lemma \ref{lem:ck-gamma-delta} that $\{\gamma_\alpha\}$ can be expressed as a polynomial in
 $\left \{c_{l_i}(\mathcal E''_i)\right \}_ {1\leq l_i\leq r,\atop 1\leq i\leq k}$ 
 where $\mathcal E''_i$ is the pull back of the, possibly twisted in the sense of \S\ref{sec:twisted}, 
 Poincar\'e bundle $\mathcal E''$ on $U \times_B \mathcal C$ via the projection $p_i \times_B \mathbbm{1}_{U} : \mathcal C^k \times_B U \rightarrow \mathcal C\times_B U$. 

\begin{theorem}\label{thm:gen-bun}
    Fix co-prime pairs $(r_i,d_i)_{i=1,\dots, N}$ and let $\mathcal E_i$ be a possibly twisted Poincar\'e bundle on $ U_{\mathcal C}(r_i,d_i)\times_B \mathcal C$. Then the total Chow group  $A_*( U_{\mathcal C}(r_1,d_1) \times_B \cdots \times_B  U_{\mathcal C}(r_N, d_N))$ is $\Q$-linearly generated by 
    \begin{equation}\label{eq6}
        \left \{ p_{U*} \left ( \prod_{i=1}^{m_1}\prod_{j=1}^{l_1} c_{d^{(1)}_j} (\mathcal E_1^{(i)}) \cdots\prod_{i=1}^{m_N}\prod_{j=1}^{l_N}c_{d^{(N)}_j} (\mathcal E_1^{(i)}) \rho^* (\Gamma)\right )  \right \}_{\Gamma\in A^*(\mathcal C^{m_1}\times_B  \cdots \times_B \mathcal C^{m_N}), d_i^{(j)}\in \mathbb Z_{\geq 0}}
    \end{equation}
    where $\mathcal E_j\ri \mathcal C\times_B U_{\mathcal C}(r_j,d_j)$ are possibly twisted Poincar\'e bundles with $\E_j^{(i)}:= q_{i,j}^*\E_j, m_i:= r_i^2(g-1)+1$ are dimensions of $U_C(r_i,d_i)$, and $q_{i,j}, \rho,p_U$ are the projections
    \begin{gather*} 
    q_{i,j}:\mathcal C^{m_1} \times_B \cdots \times_B \mathcal  C^{m_N} \times_B  U_{\mathcal C}(r_1,d_1) \times _B\cdots \times_B  U_\mathcal C(r_N, d_N)\rightarrow  \mathcal C^{m_j}\times_B U_{\mathcal C}(r_j,d_j)  \xrightarrow{i^{th}\text{ factor}} \mathcal C\times_B U_{\mathcal C}(r_j,d_j), \\
    \rho:\mathcal C^{m_1} \times_B \cdots \times_B \mathcal  C^{m_N} \times_B  U_{\mathcal C}(r_1,d_1) \times _B\cdots \times_B  U_\mathcal C(r_N, d_N)\rightarrow \mathcal C^{m_1} \times_B \cdots \times_B \mathcal  C^{m_N}, \\ p_U:  \mathcal C^{m_1} \times_B \cdots \times_B \mathcal  C^{m_N} \times_B  U_{\mathcal C}(r_1,d_1) \times _B\cdots \times_B  U_\mathcal C(r_N, d_N)\rightarrow U_{\mathcal C}(r_1,d_1) \times _B\cdots \times_B  U_\mathcal C(r_N, d_N).\end{gather*}  
\end{theorem}
\begin{proof}
    For simplicity of notation, we prove the case $N=1$. 
    The general case follows from a similar argument and the fact that $[\Delta_{ U_\mathcal C(r_1,d_1) \times_B \cdots \times_B U_\mathcal  C(r_N,d_N) /B}] = \prod_{j=1}^N \pi_j ^* [\Delta_{U_{\mathcal C}(r_j,d_j)/B}]$ 
    where $\pi_j : U_\mathcal C(r_1,d_1) \times_B \cdots \times_B U_\mathcal  C(r_N,d_N)  \rightarrow U_\mathcal C (r_j,d_j) $ is the $j$th projection.     
    We have \begin{equation*}
        \Delta_{ U_\mathcal C(r,d)/B}  =\sum_{\alpha}\gamma_\alpha \circ \delta_\alpha
    \end{equation*}where $\delta_\alpha \in A^*(U_\mathcal C(r,d)\times_B \mathcal C^m)$, $\gamma_\alpha =P_\alpha \left (\dots c_{d_j}(\mathcal E^{(i)})\dots  \right )\in A^*(\mathcal C^m\times_B  U_\mathcal C(r,d))$ is a polynomial in the Chern classes of $\mathcal E^{(i)}$, the pull-back of $\mathcal E$ on $\mathcal C\times_B  U_ \mathcal C(r,d)$ to $C^{m} \times_B  U_ \mathcal C(r,d)$ via the $i$th projection. 

    Let $z\in A^*( U_ \mathcal C(r,d))$. It follows from (\ref{diagonal as sum of composition} that \begin{equation*}
\begin{aligned}
    z & = \Delta_* (z) \\ & = \sum_{\alpha} (\gamma_\alpha \circ \delta_\alpha )_* (z) \\ & = \sum_\alpha \left (P_\alpha \left (\dots c_{d_j}(\mathcal E^{(i)})\dots  \right ) \right )_* \left( (\delta_\alpha)_*(z)\right)
\end{aligned}        
    \end{equation*}Since $(\delta_\alpha)_*(z)\in A^*(C^m)$, each term in the above sum is a $\mathbb Q$-linear combination of expressions as in (\ref{eq6}), and this completes the proof.
\end{proof}

\section{Universal curves on universal hyperelliptic Jacobians}
\subsection{Powers of the universal hyperelliptic curve}
\begin{definition}
    Let $\mC_{2,g}^n:= \times_{\h_g}^n\h_{g,1}$ be the $n$-th fibered power of the universal hyperelliptic curve over the hyperelliptic locus $\h_g\hr \m_g$. 
    Additionally, let $\psi_i:= p_i^*K\in A^1(\mC_{2,g}^n)$ be the pullback of the relative canonical bundle $K$ of $\h_{g,1}\ri \h_g$, and let $\Delta_S$ be the $S$-diagonal for any subset $S\subset [n]$. 
\end{definition}
\begin{lemma}\label{lem:gens}
    Let $\mathcal{P}$ be the set of partitions of $[n]$. For any $P\in \mathcal{P}$, let $\Delta_P:= \Pi_{I\in P}\Delta_I$ be the corresponding diagonal and $|P|$ be the number of parts of the partition. 
    For $n\leq 2g+6$, the total Chow group $A_*(\mC_{2,g}^n)$ is linearly generated by $\{\Delta_P, \psi_i\Delta_P\}_{P\in\mathcal{P},1\leq i\leq n}$.
\end{lemma}
\begin{proof}
  Let the strict diagonal of $\Delta_P$ be denoted by 
  $$\Delta_P^0:=\{(C,x_1,\cdots,x_n)|x_i=x_j \iff \exists\, S\in P \text{ such that } i,j\in S\}.$$
We can stratify $\mC_{2,g}^n$ as
$$\mC_{2,g}^n = \bigsqcup_{P\in \mathcal{P}} \Delta_P^0.$$
The strict-diagonal $\Delta_P^0$ is naturally isomorphic to $\h_{g,|P|}$. 

By \cite[Corollary 1.5]{cl-hyp}, we have that $A^*(\h_{g,m})$ is linearly generated by $[\h_{g,m}]$, and $\psi$-classes for $m\leq 2g+6$. 
The lemma now follows from a straightforward excision argument.
\end{proof}
\begin{definition}
    For $n\geq 2g+2$ and $S=\{a_1,\cdots,a_{2g+2}\}\subset [n]$, let $\alpha_S^{tav} := \underset{I_S}{\Sigma} b_{i_1,i_2}\cdots b_{i_{2g+1},i_{2g+2}}$ be the cycle defined by Tavakol \cite[(2)]{mehdi}, where $b_{ij}:= \Delta_{ij}-\dfrac{\psi_i+\psi_j}{(2g-2)}$ and $I_S$ is the set of partitions of $S$ into $g+1$-subsets with two elements each.
\end{definition}
\begin{theorem}\label{thm:c-powers}
The families $\mC_{2,g}^n\ri\h_g$ satisfy SFP for $1\leq n\leq 2g+6$. In particular, we have
  \begin{equation*} \begin{aligned}A^*(\mC_{2,g}^n) = \mathbb{Q}[\Delta_{ab},\psi_i]_{1\leq i,a,b\leq n,\atop a,b \text{ distinct}} & \bigg / \bigg \langle \Delta_{ab}(\psi_a-\psi_b),\psi_a((2g-2)\Delta_{ab}-\psi_b),
  (2g-2)\Delta_{ab}^2+\psi_a\psi_b,\Delta_{ab}-\Delta_{ba},\psi_a^2, \\ & \Delta_{ab}\Delta_{bc} - \dfrac{\psi_a\Delta_{bc}+\psi_b\Delta_{ac}+\psi_c\Delta_{ab}}{2g-2}+\dfrac{\psi_a\psi_b+\psi_b\psi_c+\psi_a\psi_c}{(2g-2)^2},\alpha_S^{tav}\bigg\rangle _{\substack{1\leq a,b,c\leq n ,\\ a,b,c \text{ distinct} , \\ S\subset [n],|S|=2g+2}}.
\end{aligned}
\end{equation*}
\end{theorem}
\begin{proof}
We begin with verifying the claimed relations. It is easy to see that 
$$\Delta_{abc} = \Delta_{ab}\Delta_{bc} = \Delta_{ab}\Delta_{bc}.$$
The claimed relations are proven in \cite[Theorem 3.10 - Corollary 3.3]{mehdi}.
For\\ $S=\{a_1,a_2,\cdots,a_t\}\subset[n]$, it is easy to see that $\Delta_S = \Pi_{2\leq i\leq t}\Delta_{a_1a_i}$.
By Lemma \ref{lem:gens}, we have the claimed generators for $A^*(\mC_{2,g}^n)$. 

Let $s:\spec \mathbb{C}\ri \h_g$ be a geometric point corresponding to a hyperelliptic curve $C$. Using \cite[Propostion 3.16 - Corollary 3.18]{mehdi}, we have that the restriction $A^*(\mC_{2,g}^n)\ri H^*(C^n)$ is injective. Hence, we have SFP for $\mC_{2,g}^n\ri \h_g$.
\end{proof}
\begin{corollary}\label{cor:c-powers}
    The total Chow group $A_*(\mC_{2,g}^n)$ is $\Q$-linearly generated by $$\{\Pi_{i\in S} \psi_i.\Pi_{1\le j\leq k} \Delta_{a_{2j-1},a_{2j}}\}_{k\geq 0,\substack{S\cup\{a_1,\cdots,a_{2k}\}\subset [n],S\cap\{a_1,\cdots,a_{2k}\}=\phi}}.$$
\end{corollary}
\begin{rem} \label{Gross-Schoen}
    The aforementioned relation:
	$$\Delta_{ab}\Delta_{bc} - \dfrac{\psi_a\Delta_{bc}+\psi_b\Delta_{ac}+\psi_c\Delta_{ab}}{2g-2}+\dfrac{\psi_a\psi_b+\psi_b\psi_c+\psi_a\psi_c}{(2g-2)^2},$$
	 is often referred to as the \textit{Gross-Schoen} cycle. For more details, see \cite{schoen,zhang-gs}.
\end{rem}
\subsection{Universal hyperelliptic Jacobians and SFP}
We study SFP for several families that arise as factors in the motivic decompositions introduced in \S \ref{sec:motive-families}, and further studied in \S\ref{sec:motivic-decomp}.
\begin{theorem}\label{thm:sfp}
        For $1\leq n\leq g+6$ and $d,d'\in\Z$, the stacks $\mC_{2,g}^n\times_{\h_g}\jac^d_{2,g}, \mC_{2,g}^n\times_{\h_g}\J^d_{2,g}, \jac^d_{2,g}\times_{\h_g} \jac^{d'}_{2,g}$, and $\J^d_{2,g}\times_{\h_g} \J^{d'}_{2,g}$ satisfy SFP over $\h_g$.
\end{theorem}
\begin{proof}
We show that $\mC_{2,g}^n\times_{\h_g}\jac^d_{2,g}$ satisies SFP. 
The claim would then follow for $\mC_{2,g}^n\times_{\h_g}\J^d_{2,g}$ by Theorem \ref{thm:gm-sfp} since $\J^d_{2,g}\ri \jac^d_{2,g}$ is a torsion-$\mathbb{G}_m$-gerbe \cite[Proposition 6.6]{melo}.

We can find integers $a_1\cdots,a_g\in \Z \setminus 0$ such that $\Sigma a_i = d$. Let $\otimes_a:\mC_{2,g}^g\ri \jac^d_{2,g}$ be the map induced by the chosen tuple $(a_1,\cdots,a_g)$. 

The induced map $1^n\times \otimes_a: \mC_{2,g}^{n+g}\ri \mC_{2,g}^n\times_{\h_g}\jac^d_{2,g}$ is generically finite and surjective. Consequently, we have $(1^n\times \otimes_a)^*: A^*(\mC_{2,g}^n\times_{\h_g}\jac^d_{2,g})\ri A^*(\mC_{2,g}^{n+g})$ is injective, and we have SFP for $\mC_{2,g}^n\times_{\h_g}\jac^d_{2,g}\ri \h_g$ by Theorem \ref{thm:c-powers} and Theorem \ref{thm:inj-sfp}.

For $\jac^d_{2,g}\times_{\h_g} \jac^{d'}_{2,g}$, we fix coordinate-wise non-zero tuples $a= (a_1,\cdots,a_g),b=(b_1,\cdots, b_g)\in\Z^g$ such that $d=\Sigma a_i, d'=\Sigma b_j$. The induced map defined by $\otimes_a\times\otimes_b: \mC^{2g}_{2,g}\ri \jac^d_{2,g}\times_{\h_g} \jac^{d'}_{2,g}$ is generically finite and surjective. Repeating a similar set of arguments as for $\mC_{2,g}^n\times_{\h_g}\jac^d_{2,g}$ above, we have SFP for $\jac^d_{2,g}\times_{\h_g} \jac^{d'}_{2,g}\ri \h_g$. A straightforward two-step application of Theorem \ref{thm:gm-sfp} gives SFP for $\J^d_{2,g}\times_{\h_g} \J^{d'}_{2,g}$.
\end{proof}

\begin{corollary}\label{cor:sfp-cpj}
	 For $1\leq n\leq g+6$ and $d\in\Z$, the stacks $[\mC_{2,g}^n/S_n]\times_{\h_g}\jac^d_{2,g}$ satisfy SFP over $\h_g$.
\end{corollary}
\begin{proof}
	The pullback along the symmetrising map $\times_{\h_g}^n\mC_{2,g}\times_{\h_g}\jac^d_{2,g}\ri [\times_{\h_g}^n\mC_{2,g}/S_n]\times_{\h_g}\jac^d_{2,g}$ 
	induces an injection on Chow groups. The claim follows from Theorem \ref{thm:inj-sfp}.
\end{proof}
\begin{rem} \label{universal Picard satisfies SFP over hyperelliptic locus}
    A similar argument as above shows that $\mathscr J^d_{2,g}\rightarrow \mathcal H_g$ also satisfies SFP.
\end{rem}
\begin{definition}\label{def:theta-xi}
    For any $d\in\Z$, let $\Xi_d$ be the Poincar\'e bundle on $\mC_{2,g}\times_{\h_g}\J^d_{2,g}$, and let $p_\C,p_\J$ be projections onto $\mC_{2,g},\J^d_{2,g}$ respectively. Additionally, let $p_\mC,p_\jac$ be projections of $\mC_{2,g}\times_{\h_g}\jac^d_{2,g}$ onto $\mC_{2,g}, \ \jac^d_{2,g}$ respectively.
\end{definition}
\begin{lemma}\label{lem:Kxi}
        The following cycle is in the image of the pullback $(1\times\nu_{1,d})^*A^1(\mC_{2,g}\times_{\h_g}\jac^d_{2,g})$:
    $$ \Xi_d-\dfrac{p_{\J}^*p_{\J *}(\Xi_d.p_\C^*K)}{2g-2}.$$
     We shall denote this cycle by $\xi_d$. Moreover, we have $\xi_d.p_\mC^*K = 0$. 
\end{lemma}
\begin{proof}
We start by showing that $\Xi_d-\dfrac{p_{\J}^*p_{\J *}(\Xi_d.p_\C^*K)}{2g-2}\in Im(1\times\nu_{1,d})^*$. For a family 
$\lb\ri\C\xrightarrow{\pi_S} S$ defining a morphism $S\ri \jac^d_{2,g}$ and $\lb_S\in Pic(S)$: 
$$\lb+\pi_S^*\lb_S-\dfrac{\pi_{S}^*\pi_{S *}((\lb+\pi_S^*\lb_S).K_{\pi_S})}{2g-2} = \lb-\dfrac{\pi_{S}^*\pi_{S *}(\lb.K_{\pi_S})}{2g-2}.$$
    We have the claimed containment in $Im(1\times\nu_{1,d})^*$ by \cite[Theorem 2.10]{shubham-r2}. 

Next, we show that $\xi_d.p_\mC^*K = 0$. We verify the relation over the cover $\h^w_{g,1}\ri \h_g$.
With slight abuse of notation, we shall identify all previously defined cycles with their respective pullbacks, 
and denote the projection $\mC_{2,g,1}^w\times_{\h_{g}}\J^d_{2,g}\ri \h^w_{g,1}\times_{\h_g}\J^d_{2,g}$ by $p_\J$.

Let $i_\s:\s_{g,1,1}\times_{\h_{g}}\J^d_{2,g}\hr \mC_{2,g,1}^w\times_{\h_g}\J^d_{2,g} $ be the inclusion map.
 It is enough to check that $$\xi_d|_{\s_{g,1,1}\times_{\h_{g}}\J^d_{2,g}} = 0.$$
 
Before we verify the triviality of the claimed restriction, we state the following relation on $\mC_{2,g,1}$ 
$$p_\mC^*K = (2g-2)\s_{g,1,1}$$ 
since $A^1(\h^w_{g,1})=0$. Consequently, we have 
\begin{equation*}i_\s^*\Xi_d = p_{\J *}i_{\s *}i_\s^*\Xi_d = p_{\J *}(\s_{g,1,1}.\Xi_d) = i_{\s}^*p_{\J}^*p_{\J *}(\s_{g,1,1}.\Xi_d)=i_\s^*\dfrac{p_{\J}^*p_{\J *}(\Xi_d.p_\C^*K)}{2g-2}.\end{equation*}
\end{proof}
\section{Universal moduli of rank $2$ bundles}\label{sec:motivic-decomp}
\noindent Following \S \ref{sec:rel-motive-decomp}, we fix $d$ odd for this section.
\subsection{Slope-stable bundles}

\begin{theorem}\label{thm:sfp-u}
	The family $\mathcal{U}^\h(2,d,g)\ri \h_g$ satisfies SFP.
\end{theorem}
\begin{proof}
	The stack $\h_g$ is smooth, Deligne-Mumford, and has a global presentation as the group quotient of a smooth scheme.
Applying Theorem \ref{thm:decomp-rel} and a standard Edidin-Graham approximation argument (see \cite[Definition-Proposition 1, and Propositions 14-15]{edidin-graham}), we get the following motivic decomposition in $\mathrm M(\mathcal H_g)$:
\begin{align*}       
\mathfrak h_{\h_g} (\mathcal{U}^\h(2,d,g))\simeq \mathfrak h_{\h_g} (\jac^d_{2,g}) \otimes  \Biggl(\mathfrak h _{\h_g}( \mathcal C^{(g-1)}_{2,g})(g-1) \oplus  \bigoplus_{i=0}^{g-2} \mathfrak h_{\h_g}(\mathcal C_{2,g}^{(i)}) \otimes 
 \biggl(\mathbb Q(i) \oplus \mathbb Q(3g-3-2i)  \biggr)\Biggr)
\end{align*}
where $\mC_{2,g}^{(n)}:= [\mC_{2,g}^n/S_n]$ for all $n\geq 0$. Consquently, we obtain a corresponding decomposition of Chow groups. The claim now follows from Corollary \ref{cor:sfp-cpj}.
\end{proof}

\begin{corollary}\label{cor:sfp-u}
	The family $\U^\h(2,d,g)\ri \h_g$ satisfies SFP.
\end{corollary}
\begin{proof}
The rigidification map $\U^\h(2,d,g)\ri \mathcal{U}^\h(2,d,g)$ is a torsion-$\gm$-gerbe \cite[Theorem B]{brauer}. 
The corollary follows from Theorem \ref{thm:sfp-u} and Theorem \ref{thm:gm-sfp}.
\end{proof}

\subsection{Filtration stacks and the excision}
We show that excision sequences on $A^*(\ve^\h(2,d,g))$ given by the Harder-Narasimhan (HN) stratification are left exact.
\begin{definition}[\cite{prada,shubham-r2-ss}]
We shall represent universal filtration stacks over $\h_g$ of type $(1,2)$ and degree $(n,d)$ as $\ext_{(1,2)}^{(n,d)}(g)$.
\end{definition}
\begin{theorem}\label{thm:excision}
	The HN-stratification of $\ve^\h(2,d,g)$ induces the following:
	$$\ve^\h(2,d,g) = \U^\h(2,d,g)\bigsqcup\left (\bigsqcup_{n\geq \lfloor d/2\rfloor +1}\ext_{(1,2)}^{(n,d)}(g) \right )$$
	such that the pushforward along the closed immersion
	$$A_*(\ext_{(1,2)}^{(n_0,d)}(g))\ri A_*\left(\ve^\h(2,d,g)\setminus \bigsqcup_{n\geq n_0+1}\ext_{(1,2)}^{(n,d)}(g)\right)$$
is injective for all $n_0\geq \lfloor d/2\rfloor +1$.
\end{theorem}
\begin{proof}
	The HN-stratification in terms of HN-polygons has been studied extensively. For more details, see \cite[\S 5]{MR3013029}.
The map 
$$\ext_{(1,2)}^{(n_0,d)}(g)\ri  \ve^\h(2,d,g)\setminus \bigsqcup_{n\geq n_0+1}\ext_{(1,2)}^{(n,d)}(g)$$ is a closed immersion by \cite[Remark 5.2]{MR3013029}.

Applying Theorem \ref{thm:sfp}, we have SFP for the families $\ext_{(1,2)}^{(n,d)}(g)$ for all $n\geq \lfloor d/2\rfloor +1$ since the forgetful morphism $\ext_{(1,2)}^{(n,d)}(g) \ri \J^n_{2,g}\times_{\h_g}\J^{d-n}_{2,g}$ is a vector bundle stack of relative dimension $g-1+d-2n$ \cite[Corollary 3.2]{prada}.
Therefore, we have the following diagram of excision sequences for any hyperelliptic curve $C$ of genus $g$:
$$\begin{tikzcd}
	A^*(\ve^\h(2,d,g)\setminus \bigsqcup_{n\geq n_0}\ext_{(1,2)}^{(n,d)}(g))\ar[r]&H^*(Bun_{2,d}^{\leq (n_0,d)})\\
	A^*(\ext_{(1,2)}^{(n_0,d)}(g))\ar[r,,hook]\ar[u,"i_{g,n_0 *}"]&H^*(\ext_{(1,2)}^{(n,d)}(C))\ar[u,"i_{C,n_0*}",hook]
\end{tikzcd}$$
We have injectivity of $i_{C,n_0*}$ by \cite[Corollary 5.14]{MR3013029}. 
Consequently, we get injectivity of $i_{g,n_0*}$.
\end{proof}
We have the following corollary as an immediate consequence of the Theorem.
\begin{corollary}\label{cor:decomp-ve}
	We have the following decomposition of the total Chow group of $\ve^\h(2,d,g)$:
	$$A_*(\ve^\h(2,d,g)) = A_*(\U^\h(2,d,g)\bigoplus \left (\bigoplus_{n\geq \lfloor d/2\rfloor +1} A_*(\ext_{(1,2)}^{(n,d)}(g)) \right ).$$
\end{corollary}
\begin{theorem}\label{thm:sfp-ve}
	The family $\ve^\h(2,d,g)\ri \h_g$ satisfies SFP.
\end{theorem}
\begin{proof}
We have SFP for $\U^\h(2,d,g)$ by Corollary \ref{cor:sfp-u}. For $\ext_{(1,2)}^{(n,d)}(g)$, SFP was shown in Theorem \ref{thm:excision}.
The claim now follows from Corollary \ref{cor:decomp-ve}.
\end{proof}

\begin{theorem}\label{thm:tautological-ve}
	Given $d\in\Z$ odd, the tautological rings $R^*(\ve^\h(2,d,g))$ and $R^*(\U^\h(2,d,g))$ are generated by the twisted-kappa classes $\{\ka_{-1,2,0},\ka_{0,1,0},\ka_{-1,0,1},\ka_{0,0,1},\ka_{-1,1,1},\ka_{-1,0,2}\}$ for all $g\geq 2$.
Additionally, we have the following relations in $A^*(\ve^\h(2,d,g))$:
\begin{align*}
   \bullet \quad& (2g-2)^{n-1}\ka_{-1,0,n} =\ka_{0,0,1}^{n-2}(n(n-1)(g-1)\ka_{-1,0,2}-(n^2-2n)\ka_{-1,0,1}\ka_{0,0,1}) \quad\forall n\geq 2,\\
  \bullet\quad & (2g-2)^{n+m-2} \ka_{-1,m,n} = {(2g-2)^{m-2}}\ka_{0,0,1}^n\ka_{-1,m,0}+{(2g-2)^{n-2}}\ka_{0,1,0}^m\ka_{-1,0,n}+\\
 & \qquad\qquad\qquad\qquad\qquad\qquad nm\ka_{0,1,0}^{m-1}\ka_{0,0,1}^{n-1}(\ka_{-1,1,1}-\dfrac{1}{2g-2}(\ka_{-1,0,1}\ka_{0,1,0}+d\ka_{0,0,1}))\quad \forall m,n\geq 1,
\\\bullet \quad &(2g-2)^{m+n-1}\ka_{0,m,n} = \ka_{0,1,0}^m\ka_{0,0,1}^n\quad \forall m,n\geq 0,\\ 
	\bullet\quad &\ka_{i,m,n} = 0 \quad \forall i\geq 1.
\end{align*}
\end{theorem}

\begin{proof}
	By excision on $\U^\h(2,d,g)\subset \ve^\h(2,d,g)$ and Theorem \ref{thm:sfp-ve}, it is enough to verify the claimed relations in 
	$H^*(Bun_{2,d}C)$ for a fixed hyperelliptic curve $C$ of genus $g$. 

Following the notation in \cite[Theorem 3.6.1]{shubham-r2-ss}, let $\E_{d,C}\ri C\times Bun_{2,d}C$ be the universal bundle and for $i\in\{1,2\}$, let
$$c_i(\E_{d,C}) := 1\otimes a_i+\Sigma_{j=1}^{2g}\gamma_j\otimes b_{j}^{i}+[pt]\otimes f_i$$
be the K\"unneth decomposition of the Chern classes of $\E_{d,C}$ where $\{\gamma_j\}_{1\leq j\leq 2g}$ is a symplectic basis of $H^1(C)$.
We have $$a_1 = \dfrac{\ka_{0,1,0}}{2g-2},a_2 = \dfrac{\ka_{0,0,1}}{2g-2}, f_1 = d, f_2 = \ka_{-1,0,1}.$$

It immediately follows that $\ka_{i,m,n} = 0$ for all $i\geq 1$ since $[pt]^2 = 0 \in H^*(C)$. 
On the other hand, we have the following relations in $H^*(Bun_{2,d}C)$:
\begin{itemize}
    \item $\ka_{0,m,n} = (2g-2)a_1^ma_2^n = (2g-2).\dfrac{\ka_{0,1,0}^m}{(2g-2)^m}.\dfrac{\ka_{0,0,1}^n}{(2g-2)^n} = \dfrac{\ka_{0,1,0}^m\ka_{0,0,1}^n}{(2g-2)^{m+n-1}}\quad\forall m,n \geq 0,$
    \item $
        \ka_{-1,0,n} = na_2^{n-1}f_2-2{n\choose 2}a_2^{n-2}\bigg(\displaystyle{\sum_{j=1}^g}b_{j}^{2}b_{j+g}^{2}\bigg) \\ 	={\dfrac{n\ka_{0,0,1}^{n-1}\ka_{-1,0,1}}{(2g-2)^{n-1}}+\frac{n(n-1)\ka_{0,0,1}^{n-2}}{2(2g-2)^{n-2}}\left(\ka_{-1,0,2}-\dfrac{\ka_{0,0,1}\ka_{-1,0,1}}{g-1}\right) \forall n\geq 2}$,
    \item $\ka_{-1,m,n} = a_1^m\ka_{-1,0,n}+a_2^n\ka_{-1,m,0} - mna_1^{m-1}a_2^{n-1} \bigg(\displaystyle\sum_{j=1}^gb_{j}^{1}b_{j+g}^{2}-b_{j+g}^{1}b_{j}^{2}\bigg )  \\
=\dfrac{\ka_{0,1,0}^m\ka_{-1,0,n}}{(2g-2)^m}+
\dfrac{\ka_{0,0,1}^n\ka_{-1,m,0}}{(2g-2)^n}	+\dfrac{nm\ka_{0,1,0}^{m-1}\ka_{0,0,1}^{n-1}}{(2g-2)^{m+n-2}}\left( \ka_{-1,1,1}-\dfrac{\ka_{0,1,0}\ka_{-1,0,1}+\ka_{0,0,1}\ka_{-1,1,0}}{2g-2} \right) \\ \forall m,n\geq 1$.
\end{itemize}
Applying \cite[Corollary 1.2]{larson-24}, we have the claimed generators for $R^*(\ve^\h(2,d,g))$.
\end{proof}

\noindent We end this section with a remark which we shall refer to repeatedly in the rest of the paper.
\begin{rem}\label{monodromy argument}
    Fix $C$, a hyperlliptic curve in $\mathcal H_g$. $\mathfrak X\rightarrow \mathcal H_g$ is smooth where $\mathfrak X$ is one of the stacks considered in Theorems \ref{thm:sfp} , \ref{thm:sfp-u} and \ref{thm:sfp-ve}. Global cycles in all of these cases give monodromy invariant cohomology after restricting to a fiber over a geometric point. 

    \noindent We elaborate more on this point. Let $\mathfrak X \xrightarrow{\pi} B$ be a smooth morphism between smooth algebraic stacks and $B$ is Deligne-Mumford. Then the cycle class map gives a group homomorphism $cl: A^r (\mathfrak X) \rightarrow H^{2r}(\mathfrak X,\mathbb Q) \simeq \operatorname{Hom}_{D(B)} \big ( \mathbb Q_B, R \pi_* \mathbb Q_{\mathfrak X} [2r] \big )$. Taking $0$th cohomology gives us a morphism $A^r (\mathfrak X) \rightarrow H^0(B, R^{2r}\pi_* \mathbb Q_{\mathfrak X})$. Suppose additionally, that $R^{2r}\pi_* \mathbb Q_\mathfrak X$ is a local system with stalk $H^{2r}(\mathfrak X_b, \mathbb Q)$. Then we have the isomorphism $H^0 (B, R^{2r}\pi_* \mathbb Q_\mathfrak X) \simeq H^{2r} (\mathfrak X_b, \mathbb Q)^{\pi_1(B,b)}$, so we finally get a homomorphism \begin{equation*}
        A^{r}(\mathfrak X) \rightarrow H^{2r}(\mathfrak X_b,\mathbb Q )^{\pi_1(B,b)}.
    \end{equation*}  
    In particular if $\mathfrak X\xrightarrow{\pi}B$ is topologically locally trivial upto homotopy type, the cycle class map and restriction to a fiber gives a homomorphism $ A^{r}(\mathfrak X) \rightarrow H^{2r}(\mathfrak X_b,\mathbb Q )^{\pi_1(B,b)}$.

    \noindent For smooth proper morphisms this follows from the Leray spectral sequence and Ehresmann's theorem, and for torsion $\mathbb G_m$-gerbes over a smooth proper family this can be shown using Leray-Hirsch. Finally for $\mathcal Vec^\mathcal H(2,d,g) \rightarrow \mathcal H_g$, this follows from the works of Atiyah-Bott in \cite{atiyah-bott}. So, in view of SFP, the cycle class map and restriction to a geometric fiber gives an embedding $$A^* (\mathfrak X) \hookrightarrow H^{2*} (\mathfrak X_{[C]},\mathbb Q)^{\pi_1(\mathcal H_g, [C])}.$$  We know that the standard monodromy representation on $H^{*}(\mathfrak X_{[C]},\mathbb Q)$ in each of these cases is induced by an algebraic representation of $\operatorname{Sp}(V)$ where $V=H^1(C,\mathbb Q)$ with the standard action. Moreover, from the results of \cite{MR535062} (see also \cite[Theorem 8.1]{MR4182796}), we know that the image of $\pi_1(\mathcal H_g ,[C])$ under the standard monodromy representation is Zariski dense in $\operatorname{Sp}(V)$ . This gives us $$H^{2*} (\mathfrak X_{[C]},\mathbb Q)^{\pi_1(\mathcal H_g, [C])}=H^{2*} (\mathfrak X_{[C]},\mathbb Q)^{\operatorname{Sp}(V)}.$$
\end{rem}
\section{Proofs of the main theorems}\label{sec:conj-hyp}
\subsection{Product of universal Jacobians}
The following results are analogues of \cite[Theorem 4.4]{shubham-r2}. We include them below for completeness and to address minor technical difficulties that arise when we generalize the result to higher genus.

We also indicate here how this gives an alternative proof of \cite[Theorem 1.1]{larson-24}. We use the fact for any $d\in\Z$, there exists $\T_d\in A^1(\jac^d_{2,g})$ such that the restriction of $\T_d$ to $J^d_C\simeq J^{g-1}_C$ is $\T_C$ in $H^2(J^{g-1}_C,\Q)$ where the isomorphism $J^d_C\simeq J^{g-1}_C$ is given by twisting with an arbitrary divisor of degree $d-g+1$. This follows from \cite[Corollary 6.3(ii)]{melo}. Moreover, we can choose $\Theta_d$ so that
\begin{equation} \label{defintion of Theta_d}\T_d := \dfrac{d\ka_{0,1}-(g-1)\ka_{-1,2}}{2(g-1)}\in A^1(\jac^d_{2,g}).\end{equation}
\begin{theorem}\cite[Theorem 1.1]{larson-24} \label{simpler proof of Larson}
    We have the following isomorphism of graded rings \begin{equation*}
        A^* (\mathscr J^d_{2,g}) \simeq H^{2*} (\mathscr Pic^d(C),\mathbb Q)^{\pi_1(\mathcal H_g,[C])} \simeq  \mathbb Q[\kappa_{0,1},\kappa_{-1,2}] / \langle (d\kappa_{0,1} - (g-1)\kappa_{-1,2})^{g+1} \rangle 
    \end{equation*}
\end{theorem}
\begin{proof}
    We know $\mathscr J^{d}_{2,g} \rightarrow \mathcal J^d_{2,g}$ is a torsion $\mathbb G_m$-gerbe so by SFP and Remark \ref{monodromy argument}, we already get an embedding of graded rings $A^*(\mathscr J^d_{2,g}) \hookrightarrow H^{2*}(\mathscr Pic^d(C), \mathbb Q)^{\pi_1(\mathcal H_g, [C])}=H^{2*}(\mathscr Pic^d(C), \mathbb Q)^{\operatorname{Sp}(V)} $ where $V=H^1(C,\mathbb Q)$. So we only need to show surjectivity and compute the ring of invariants. From the works of \cite{MR3013029}, we know that there is a graded ring isomorphism \begin{equation*}
        H^* (\mathscr Pic^d(C),\mathbb Q) \simeq \bigwedge^* V \otimes \mathbb Q[u], \ \deg u=2
    \end{equation*} which is also $\operatorname{Sp}(V)$-equivariant (the action of $\operatorname{Sp}(V)$ on $\mathbb Q[u]$ is trivial). So we have \begin{equation*}
        H^{*}(\mathscr Pic^d(C), \mathbb Q)^{\operatorname{Sp}(V)} = \mathbb Q[\theta]/(\theta^{g+1}) \otimes \mathbb Q[u]
    \end{equation*}    where $\theta \in \bigwedge^2 V$ is dual to the symplectic form. In fact, from \cite[Theorem 4.7]{MR3013029} and the calculations in loc. cit., we have the formulas $\theta = \dfrac{d\kappa_{0,1}-(g-1)\kappa_{-1,2}}{2(g-1)}$ and $u=\dfrac{\kappa_{0,1}}{2(g-1)}$. This completes the proof of surjectivity since the $\kappa$-classes are global cycles on $\mathscr J^d_{2,g}$. 
\end{proof}
\begin{lemma}\label{lem:jxj-gen}
    For all $g\geq 2$, the ring $A^*(\jac^0_{2,g}\times_{\h_g}\jac^{0}_{2,g})$ is generated by $\pi_1^*\T_0,\pi_2^*\T_0,\mu^*\T_0$ as a $\Q$-algebra where $\mu:\jac^0_{2,g}\times_{\h_g}\jac^{0}_{2,g}\ri \jac^0_{2,g}$ is the twisting map and $\pi_i$ is the projection onto the $i$-th factor for $i=1,2$.
\end{lemma}
\begin{proof} 
Following the $\Q$-linear given by Theorem \ref{thm:gen-bun}
for tuples $\{(1,0),(1,0)\}$, along with the notation in Lemma \ref{lem:Kxi}, $A^*(\jac^0_{2,g}\times_{\h_g}\jac^{0}_{2,g})$ is $\Q$-linearly generated by the following cycles:
\begin{equation}\label{eq:gens}\big \{p_{J\times J *}\left(\Pi_{1\leq i\leq 2g}c_1(\xi_0)^{a_i}\rho^*(\Gamma)\right)\big \}_{\Gamma\in A^*(\mC_{2,g}^{2g}), a_i\geq 0}\end{equation}
where $p_{J\times J}, \rho$ are projections of $\jac^0_{2,g}\times_{\h_g}\jac^{0}_{2,g}\times_{\h_g}\mC_{2,g}^{2g}$ onto $ \jac^0_{2,g}\times_{\h_g}\jac^{0}_{2,g},\mC_{2,g}^{2g}$ respectively.

We analyze the restrictions of these generators, in (\ref{eq:gens}), over $[C]\in |\h_g|$ since it is enough to show the claim for the image of the restriction map $A^*(\jac^0_{2,g}\times_{\h_g}\jac^{0}_{2,g})\ri H^*(J^0_C\times J^0_C)$ for a hyperelliptic curve $C$ of genus $g$  by Theorem \ref{thm:sfp}. 

Fixing a symplectic basis of $H^1(C,\Q)$ and the K\"unneth formula, an immediate consequence is that $A^*(\jac^0_{2,g}\times_{\h_g}\jac^{0}_{2,g})$ is a subalgebra of 
$$\bigwedge_{\substack{1\leq i\leq 2g,\\1\leq j\leq 2}} [b^j_i].$$

Using the image of the restriction $A^*(\mC_{2,g}^{2g}) \ri H^*(C^{2g})$, given by Theorem \ref{thm:c-powers}, 
$A^*(\jac^0_{2,g}\times_{\h_g}\jac^{0}_{2,g})$ is $\Q$-linearly generated by
$$\{(\sum_{1\leq i\leq g} b^1_ib^1_{i+g})^a(\sum_{1\leq j\leq g} b^2_jb^2_{j+g})^b(\sum_{1\leq k\leq g} b^1_kb^2_{k+g}+b^2_kb^1_{k+g})^c\}_{a,b,c\geq 0}.$$
The claim follows consequently since $\pi_1^*\T_0,\pi_2^*\T_0,\mu^*\T_0$ can be identified with 
$$\sum_{1\leq i\leq g} b^1_ib^1_{i+g}, \ \sum_{1\leq j\leq g} b^2_jb^2_{j+g}, \  \sum_{1\leq k\leq g} (b^1_k+b^2_k)(b^1_{k+g}+b^2_{k+g})$$ respectively.
\end{proof}

\begin{theorem}\label{thm:jxj-rel}
 Recall the definition of $\Theta_d$ in Equation \ref{defintion of Theta_d}.
 $$x:= \pi_1^*\T_d, \ y:= \pi_2^*\T_{d'}, \ z:= \mu_{d,d'}^*\T_{d+d'}-x-y \in A^*(\jac^d_{2,g}\times_{\h_g}\jac^{d'}_{2,g})$$
  with $\pi_1,\pi_2$ projections onto the first and second factors respectively, and $\mu_{d,d'}:\jac^d_{2,g}\times_{\h_g}\jac^{d'}_{2,g} \to \jac^{d+d'}_{2,g}$ the twisting map.
 The ring $A^*(\jac^d_{2,g}\times_{\h_g}\jac^{d'}_{2,g})$ is isomorphic to
$$\Q[x,y,z]/I$$
where the ideal $I$ is generated by 
$$ \left\{{\sum}_{\substack{m+n+p=g+1,\\p+2n=t}}{{g+1}\choose{m,n,p}}x^my^nz^p\right\}_{0\leq t\leq 2g+2}.$$
\end{theorem}
\begin{proof}
A straightforward adaptation of \cite[Theorem A.1]{bae-maulik} shows that we have an isomorphism 
$$\mathfrak{h}({\jac^d_{2,g}\times_{\h_g}\jac^{d'}_{2,g}})\simeq \mathfrak{h}(\jac^0_{2,g}\times_{\h_g}\jac^{0}_{2,g})$$ preserving the cup product for all $d,d'\in \Z$. 
Therefore, it is enough to show the claim for all $d,d'=0$.
Consequently, we have that $A^*(\jac^d_{2,g}\times_{\h_g}\jac^{d'}_{2,g})$ is generated by $x,y,z$ as a $\Q$-algebra by Lemma \ref{lem:jxj-gen}.

Following the notation established in Lemma \ref{lem:jxj-gen}, set $u_N := \sum_{1\leq l\leq g}(b^1_l+Nb^2_l)(b^1_{l+g}+Nb^2_{l+g})$ for all $N\in \Z$.
It is clear that $$u_N = x+Nz+N^2y,\quad u_N^{g+1}=0\quad \forall N\in\Z.$$
The claimed generators are coefficients of $N^t$ in the expansion of $(x+Nz+N^2y)^{g+1}$ for all $0\leq t\leq 2g+2$.
Therefore, the generators are contained in $I$.

We show that $I$ is the ideal of relations for $A^*(\jac^d_{2,g}\times_{\h_g}\jac^{d'}_{2,g})$.
It is clear that $A^*(\jac^d_{2,g}\times_{\h_g}\jac^{d'}_{2,g}), \ \Q[x,y,z]/I$ are graded rings with top degree $2g$ such that $x,y,z$ are algebraically independent up to degree $g$ in both.

We start by showing Poincar\'e duality for $A^*(\jac^d_{2,g}\times_{\h_g}\jac^{d'}_{2,g})$. Set 
$$R:= \underset{\substack{1\leq i\leq 2g,\\1\leq j\leq 2}}{\bigwedge} [b^j_i].$$
 As explained in the proof of Lemma \ref{lem:jxj-gen}, 
we can identify $x,y,z$ with $$\sum_{1\leq i\leq g} b^1_ib^1_{i+g}, \sum_{1\leq j\leq g} b^2_jb^2_{j+g}, \sum_{1\leq k\leq g} b^1_kb^2_{k+g}+b^2_kb^1_{k+g}$$ respectively.
Let $V$ be a $2g$-dimensional symplectic vector space with symplectic basis $e_1,f_1,\dots, e_g,f_g$. Consider $\tilde R:= \bigwedge^\bullet (V\oplus V)$. $\operatorname{Sp}(V)$ acts diagonally on $V\oplus V$ and hence on $\tilde R$. Then we have a canonical isomorphism $R\simeq \tilde R$ given by $b^j_i \mapsto e^{(j)}_{i}, \ 1\leq i\leq g $, $b^j_i =f^{(j)}_i, \ g+1\leq i\leq 2g $. From the standard facts about the representation theory of $\operatorname{Sp}(V)$, we know $\tilde R^{\operatorname{Sp}(V)} = \mathbb Q[ x,y,z ]$ (this follows from the results in \cite[Chapter VI]{MR1488158} or \cite[Theorem 2]{MR986027}). Using Theorem \ref{thm:sfp} and Lemma \ref{lem:jxj-gen}, we have
$$A^* (\mathcal J^d_{2,g} \times_{\mathcal H_g} \mathcal J^{d'}_{2,g}) \simeq \tilde R^{\operatorname{Sp}(V)}.$$ Using that $\operatorname{Sp}(V)$ is reductive and the non-degenerate pairing on $\tilde R$ is $\operatorname{Sp}(V)$-invariant so the pairing on $\tilde R^{\operatorname{Sp}(V)}$ is also non-degenerate. This follows from the existence of the Reynold's operator $\tilde R\rightarrow \tilde R^{\operatorname{Sp}(V)}$, as shown in \cite{MR719371} (see Definition 1.5 and the subsequent discussion). This proves Poincar\'e duality on $A^* (\mathcal J^d_{2,g} \times_{\mathcal H_g} \mathcal J^{d'}_{2,g})$. 

Since $x,y, z$ are algebraically independent upto degree $g$, we have that \begin{equation*}
    \dim_{\mathbb Q} A^k (\mathcal J^d_{2,g} \times_{\mathcal H_g} \mathcal J^{d'}_{2,g}) = \binom{k+2}{2}, \ 0\leq k\leq g.
\end{equation*}Poincar\'e duality then gives us \begin{equation*}
    \dim_{\mathbb Q}A^k (\mathcal J^d_{2,g} \times_{\mathcal H_g} \mathcal J^{d'}_{2,g})= A^{2g-k} (\mathcal J^d_{2,g} \times_{\mathcal H_g} \mathcal J^{d'}_{2,g}) , \ 0\leq k\leq g. 
\end{equation*}So we get \begin{equation*}
    \dim_{\mathbb Q} A^k (\mathcal J^d_{2,g} \times_{\mathcal H_g} \mathcal J^{d'}_{2,g}) = \binom{\min{\{k,2g-k\}+2}}{2}. 
\end{equation*}
In particular, from a straightforward computation, we get $$\dim_{\mathbb Q}A^*(\mathcal J^d_{2,g} \times_{\mathcal H_g} \mathcal J^{d'}_{2,g}) = \frac{(g+1)(g+2)(2g+3)}{6}$$
We know that there is a surjection $\mathbb Q[x,y,z]/I \twoheadrightarrow A^* (\mathcal J^d_{2,g} \times_{\mathcal H_g} \mathcal J^{d'}_{2,g})$. So we only need to prove a dimension bound for isomorphism. We have the following relations in $\mathbb Q[x,y,z]/I$ : \begin{equation}\label{AI assist} \begin{aligned}
    0&= \binom{g+1}{c} x^{g+1-c}z^c  +\sum_{j=1}^{\lfloor \frac{c}{2} \rfloor} \frac{(g+1)!}{(g+1-c+j)!(c-2j)! j!} x^{g+1-c+j}y^j z^{c-2j} \end{aligned} \end{equation} \begin{equation*} \begin{aligned}0& = \binom{g+1}{c} y^{g+1-c}z^c +\sum_{j=1}^{\lfloor \frac{c}{2} \rfloor} \frac{(g+1)!}{j!(c-2j)!(g+1-c+j)!} x^j y^{g+1-c+j}z^{c-2j}.
    \end{aligned}
\end{equation*} So $\mathbb Q[x,y,z]/I $ is spanned by $\{x^ay^bz^c\mid 0\leq c\leq g , \  0\leq a,b\leq g- c  \}$. 

\noindent This gives us the upper bound \begin{equation*}
     \begin{aligned}
         \dim_{\mathbb Q} \mathbb Q[x,y,z]/I & \leq \sum_{c=0}^g (g-c+1)^2 \\ & =  \frac{(g+1)(g+2)(2g+3)}{6}.
     \end{aligned}
 \end{equation*}It follows that $\mathbb Q[x,y,z]/I \twoheadrightarrow A^*(\mathcal J^d_{2,g} \times_{\mathcal H_g} \mathcal J^{d'}_{2,g})$ is an isomorphism. In particular $\mathbb Q[x,y,z]/I$ also satisfies Poincar\'e duality. 
\end{proof}
We rephrase the above theorem for later use in the following corollary.
\begin{corollary}\label{cor:ext-alg}
    The subalgebra of $\underset{\substack{1\leq i\leq 2g,\\1\leq j\leq 2}}{\bigwedge} [b^j_i]$ generated by the elements
    $$x:= \sum_{1\leq i\leq g} b^1_ib^1_{i+g}, \quad y:= \sum_{1\leq j\leq g} b^2_jb^2_{j+g}, \quad z:= \sum_{1\leq k\leq g} b^1_kb^2_{k+g}+b^2_kb^1_{k+g}$$
    is given by $\Q[x,y,z]/I$ where the ideal $I$ is generated by
    $$ \left\{\sum_{\substack{m+n+p=g+1\\p+2n=t}}{{g+1}\choose{m,n,p}}x^my^nz^p\right\}_{0\leq t\leq 2g+2}.$$
\end{corollary}
\begin{corollary} \label{kappa computation}
  For all $g\geq 2 $ and $d$ odd, the ring $R^*(\ve^\h(2,d,g))$ is isomorphic to:
    $$\Q[\ka_{0,1,0},\ka_{0,0,1},\ka_{-1,0,1}]\otimes_\Q \Q[x,y,z]/I$$
    where $I$ is the ideal generated by $\left\{{\sum}_{\substack{m+n+p=g+1,\\p+2n=t}}{{g+1}\choose{m,n,p}}x^my^nz^p\right\}_{0\leq t\leq 2g+2}$
    and \\
    $x:= \dfrac{d\ka_{0,1,0}-(g-1)\ka_{-1,2,0}}{2(g-1)},y:= \dfrac{\ka_{-1,0,1}\ka_{0,0,1}-(g-1)\ka_{-1,0,2} }{2(g-1) },z:= \dfrac{ d\ka_{0,0,1}+\ka_{0,1,0}\ka_{-1,0,1}-(2g-2)\ka_{-1,1,1}}{2(g-1)}$.
\end{corollary}
\begin{proof}
    In continuation of Theorem \ref{thm:tautological-ve}, recall that $\E_{d,C}$ is the Poincar\'e bundle on $C\times Bun_{2,d}C$ and
\begin{gather*}
c_i(\E_{d,C}) := 1\otimes a_i+\Sigma_{j=1}^{2g}\gamma_j\otimes b_{j}^i+[pt]\otimes f_i \text{ for }i=1,2\text{ where }
\end{gather*} $$a_1 = \dfrac{\ka_{0,1,0}}{2g-2}, \ a_2 = \dfrac{\ka_{0,0,1}}{2g-2}, \ f_1 = d, \ f_2 = \ka_{-1,0,1}.$$
Consequently, we have
\begin{align*}
    \ka_{-1,2,0} = &\dfrac{d\ka_{0,1,0}}{g-1}-2 \left (\sum_{1\leq i\leq g}b^1_ib^1_{i+g} \right ),\\
    \ka_{-1,1,1} =& \dfrac{\ka_{0,1,0}\ka_{-1,0,1}}{2g-2}+\dfrac{d\ka_{0,0,1}}{2g-2}-\left (\sum_{1\leq j\leq g}b^1_jb^2_{j+g}+b^2_jb^1_{j+g} \right),\\
    \ka_{-1,0,2} = &\dfrac{\ka_{-1,0,1}\ka_{0,0,1}}{g-1}-2 \left (\sum_{1\leq k\leq g}b^2_kb^2_{k+g} \right ).
\end{align*}
Thus, we get the claimed expressions for $x,y,z$ based on definitions in Corollary \ref{cor:ext-alg}.
\end{proof}
\begin{theorem}\label{thm:conj-odd}
    Conjecture \ref{conj:hyp} holds for all $g\geq 2$ and $d$ odd.
\end{theorem}
\begin{proof}
    By Remark \ref{monodromy argument} and SFP for $\mathcal Vec^\mathcal H(2,d,g) \rightarrow \mathcal H_g$, we have the embedding \begin{equation*}
        A^*(\mathcal Vec^\mathcal H(2,d,g)) \hookrightarrow H^{2*} (Bun_{2,d}(C),\mathbb Q)^{\operatorname{Sp}(V)},
    \end{equation*}so we only need to show that the above map is surjective. By the results of \cite{atiyah-bott} (see also \cite[Theorem 4.7]{MR3013029}), we have \begin{equation*}
        H^{2*}(Bun_{2,d}(C),\mathbb Q) =\mathbb Q[a_1,a_2,f_2] \otimes \left ( \bigwedge^\bullet (V\oplus V) \right ) 
    \end{equation*}    as $\operatorname{Sp}(V)$ representations where $\operatorname{Sp}(V)$ acts trivially on $\mathbb Q[a_1,a_2, f_2]$. By a similar argument as in Theorem \ref{thm:jxj-rel}, we have \begin{equation*}
        H^{2*} (Bun_{2,d}(C),\mathbb Q)^{\operatorname{Sp}(V)} = \mathbb Q[a_1,a_2,f_2] \otimes \mathbb Q[x,y,z]/I.
    \end{equation*} The result now follows from the computations in Theorem \ref{thm:tautological-ve} and Corollary \ref{kappa computation}. In fact this shows \begin{equation*}
        A^*(\mathcal Vec^\mathcal H(2,d,g)) =R^* (\mathcal Vec^\mathcal H(2,d,g)) \simeq H^{*}_\kappa (Bun_{2,d}(C),\mathbb Q) =H^{2*}(Bun_{2,d}(C),\mathbb Q)^{\operatorname{Sp}(V)}.
    \end{equation*}We note that this also gives a proof Theorem \ref{global Chow = monodromy invariant of a fiber}.
\end{proof}

\noindent The rest of this section is devoted to computing Hilbert series of the various stacks we consider.

\subsection{Slope-unstable locus}
\begin{theorem}\label{thm:hilb-unstable}
    For $d$ odd, the Hilbert series of $A^*\left(\bigsqcup_{n\geq \lfloor d/2\rfloor+1}\ext^{(n,d)}_{(1,2)}(g) \right )$ is given by 
    $$\dfrac{{g+2\choose 2}T^g+\sum_{1\leq k\leq g}{g-k+2\choose 2}(T^{g-k}+T^{g+k})}{(1-T)^2(1-T^2)}.$$ Moreover, this is the slope-unstable locus in $\ve^\h(2,d,g)$ and has codimension $g$.
\end{theorem}
\begin{proof}
As shown in Theorem \ref{thm:excision}, the slope-unstable locus in $\ve^\h(2,d,g)$ is given by $\bigsqcup_{n\geq \lfloor d/2\rfloor+1}\ext^{(n,d)}_{(1,2)}(g)$
and the HN-stratification induces injective pushforwards on Chow groups. 
Therefore, it is enough to compute the Hilbert series of $A^*(\J^n_{2,g}\times_{\h_g}\J^{d-n}_{2,g})$ for all $n\geq \lfloor d/2\rfloor+1$ by \cite[Corollary 3.2]{prada}.    

Applying Poincar\'e duality and algebraic independence in degree at most $g$ in Theorem \ref{thm:jxj-rel}, 
the Hilbert series of $A^*(\jac^n_{2,g}\times_{\h_g}\jac^{d-n}_{2,g})$ is given by 
$${g+2\choose 2}T^g+\sum_{1\leq k\leq g}{g-k+2\choose 2}(T^{g-k}+T^{g+k}).$$
Since $\J^n_{2,g}\times_{\h_g}\J^{d-n}_{2,g}$ can be viewed as a tower of torsion-$\gm$-gerbes over $\jac^n_{2,g}\times_{\h_g}\jac^{d-n}_{2,g}$, the Hilbert series of $A^*(\J^n_{2,g}\times_{\h_g}\J^{d-n}_{2,g})$ is given by
    $$\dfrac{{g+2\choose 2}T^g+\sum_{1\leq k\leq g}{g-k+2\choose 2}(T^{g-k}+T^{g+k})}{(1-T)^2}.$$ 
The claim now follows from the codimension of the filtration strata $\ext^{(n,d)}_{(1,2)}(g)$ in $\ve^\h(2,d,g)$ 
being $g-1+2n-d$ for all $n\geq \lfloor d/2\rfloor+1$.
\end{proof}

\subsection{The $\ka$-cohomology}
We compute the Hilbert series of $H^*_\ka(Bun_{2,d}C)$ where $C$ is a curve of genus $g$.
\begin{theorem}\label{thm:ext-hilb}
    For the graded exterior algebra 
    $$\bigwedge_{\substack{1\leq j\leq 2g,\\1\leq i\leq 2}} [b^i_j]$$
with $\deg (b^i_j)=2i-1$ for all $1\leq i\leq 2$, $1\leq j\leq 2g$,
    the Hilbert series of the subalgebra $S$ generated by
    $$\left \{x:=\sum_{1\leq i\leq g}b^1_ib^1_{i+g},\;y:=\sum_{1\leq j\leq g}b^2_jb^2_{j+g},\;z:=\sum_{1\leq k\leq g}b^1_kb^2_{k+g}+b^2_kb^1_{k+g}\right \}$$
is given by
$$\dfrac{(1-T^{g+1})(1-T^{g+2})(1-T^{2g+3})}{(1-T)(1-T^2)(1-T^3)}$$
where we half the degree of each monomial in $x,y,z$.
\end{theorem}
\begin{proof}
From the proof of Theorem \ref{thm:jxj-rel}, we already know that $S=\mathbb Q[x,y,z] \subset\displaystyle{\bigwedge_{\substack{1\leq j\leq 2g,\\1\leq i\leq 2}} [b^i_j]}$ has a basis generated by $\{x^ay^b z^c : 0\leq c\leq g , \ 0\leq a,b \leq g-c  \}$ because the ring is the same and only the grading is different. So we have a homogeneous basis of the ring. Thus we can explicitly compute \begin{equation*}
    \begin{aligned}
        \operatorname{hilb}_{\operatorname{ser}} (S) & = \sum_{c=0}^g \sum_{a,b=0}^{g-c} T^{a+3b+2c} \\ & = \frac{(1-T^{g+1})(1-T^{g+2})(1-T^{2g+3})}{(1-T)(1-T^2)(1-T^3)}
    \end{aligned}
\end{equation*} where the second equality is a nice exercise in combinatorics. Another way to see it would be to check that both satisfy the recursion $Q_g(T) =T^2Q_{g-1}(T)+ \frac{(1-T^{g+1})(1-T^{3g+3})}{(1-T)(1-T^3)}$ with $Q_0(T)=1$.  
\end{proof}

\begin{theorem}\label{thm:hilb-ka}
   The Hilbert series of $H^*_\ka(Bun_{2,d}(C),\mathbb Q)$, under half of the grading induced by $H^*(Bun_{2,d}(C),\Q)$, is given by
   $$\dfrac{(1-T^{g+1})(1-T^{g+2})(1-T^{2g+3})}{(1-T)^3(1-T^2)^2(1-T^3)}.$$ 
\end{theorem}
\begin{proof}
The generators of the ring $H^*_\ka(Bun_{2,d}C)$ were found in Theorem \ref{thm:tautological-ve}.
Following Corollary \ref{cor:ext-alg}, we have
$$H^*_\ka(Bun_{2,d}C)\simeq \Q[\ka_{0,1,0},\ka_{0,0,1},\ka_{-1,0,1}]\otimes_\Q S.$$
The degrees of $\ka_{0,1,0},\ka_{0,0,1},\ka_{-1,0,1}$ are $1,2,1$ respectively. Therefore, the Hilbert series of $H^*_\ka(Bun_{2,d}C)$ is given by
$$\hilb (H^*_\ka(Bun_{2,d}C)) = \dfrac{\hilb (S)}{(1-T)^2(1-T^2)} = \dfrac{(1-T^{g+1})(1-T^{g+2})(1-T^{2g+3})}{(1-T)^3(1-T^2)^2(1-T^3)}.$$
\end{proof}

\subsection{Slope-stable bundles}
We compute the Hilbert series of $\U^\h(2,d,g)$ for $d$ odd. 
\begin{definition}
    Let $P_{g,n}(T)$ be the Hilbert series of $\mC_{2,g}^{(n)}\times_{\h_g}\jac^d_{2,g}$ for all $n\geq 0$ where $\mC_{2,g}^{(n)}:= [\mC_{2,g}^n/S_n]$.
\end{definition}
\begin{lemma}\label{lem:stable-hilb}
    The Hilbert series of $\U^\h(2,d,g)$ is given by
    $$\dfrac{P_{g,g-1}(T)T^{g-1}+\sum_{i=0}^{g-2}P_{g,i}(T)(T^i+T^{3g-3-2i})}{1-T}.$$
\end{lemma}
\begin{proof}
The map $\U^\h(2,d,g)\ri \mcalu^\h(2,d,g)$ is a torsion-$\gm$-gerbe by \cite[Theorem B]{brauer}.
Applying Theorem \ref{thm:gerbe-ring}, we have $$\hilb(\U^\h(2,d,g)) =\dfrac{\hilb(\mcalu^\h(2,d,g))}{1-T}.$$
Plugging the motivic decomposition of $\mcalu^\h(2,d,g)$ in Theorem \ref{thm:decomp-rel}  into this expression gives the claimed equality.
\end{proof}
We now compute $P_{g,n}(T)$ for $0\leq n< g$.
\begin{notation}
The Hilbert series $P_{g,n}(T)$ is independent of the value of $d$ by an adaptation of \cite[Proposition 0.1]{bae-maulik}.
We shall refer to $P_{g,n}(T)$ as the Hilbert series of $\mC_{2,g}^{(n)}\times_{\h_g}\jac^g_{2,g}$ from here on.
\end{notation}

\begin{theorem}\label{thm:poincare-cnj}
    For $0\leq n\leq g+6$, the ring $A^*(\mC_{2,g}^{(n)}\times_{\h_g}\jac^g_{2,g})$ is of degree $n+g$ and satisfies Poincar\'e duality.
\end{theorem}
\begin{proof}

First of all let $V=H^1(C,\mathbb Q)$ which is a symplectic vector space. Then \begin{equation*}
    H^* (C^n \times J_C ,\mathbb Q) = \left ( \bigotimes_{i=1}^n \left ( \mathbb Q \oplus V \oplus \mathbb Q{\eta_i} \right ) \right )\otimes \left(\bigwedge^\bullet V \right)
\end{equation*} as $\operatorname{Sp}(V)$-modules where $\operatorname{Sp}(V)$ acts trivially on $\mathbb Q$ and $\mathbb Q\eta_i$. It follows that we have a surjective $\operatorname{Sp}(V)$-equivariant ring homomorphism \begin{equation*}
    \bigwedge^\bullet \left ( V_1\oplus \cdots \oplus V_n \oplus V_{J_C} \right ) \twoheadrightarrow  H^* (C^n \times J_C ,\mathbb Q)
 \end{equation*} and since $\operatorname{Sp}(V)$ is reductive, this induces a surjection onto their invariant spaces \begin{equation*}
     \left [\bigwedge^\bullet \left ( V_1\oplus \cdots \oplus V_n \oplus V_{J_C} \right ) \right ]^{\operatorname{Sp}(V)} \twoheadrightarrow  H^* (C^n \times J_C ,\mathbb Q)^{\operatorname{Sp}(V)}.
 \end{equation*}It follows from the standard representation theory of $\operatorname{Sp}(V)$ (see as before \cite{MR1488158} or \cite{MR986027}) that $ H^* (C^n \times J_C ,\mathbb Q)^{\operatorname{Sp}(V)}$ is generated as a ring by $ \big \{\psi_i, \Delta_{ij} ,\xi_{i,g}, p_{J_C}^*\Theta_g\big\}_{i\neq j} $. So we have established using SFP (see Theorem \ref{thm:sfp} and Remark \ref{monodromy argument}) that as graded rings \begin{equation*}
     A^* (\mC_{2,g}^{n}\times_{\h_g}\jac^g_{2,g}) \simeq  H^{2*} (C^n \times J_C ,\mathbb Q)^{\operatorname{Sp}(V)}.
 \end{equation*}
 Further taking $S_n$ -invariants, we get 
\begin{equation*}
      A^* (\mC_{2,g}^{(n)}\times_{\h_g}\jac^g_{2,g}) \simeq  H^{2*} (C^{(n)} \times J_C ,\mathbb Q)^{\operatorname{Sp}(V)}.
\end{equation*}
But the Poincare pairing on $ H^* (C^{(n)} \times J_C ,\mathbb Q)$ is non-degenerate and $\operatorname{Sp}(V)$-invariant. Since $\operatorname{Sp}(V)$ is reductive, it follows that the pairing on $ H^* (C^{(n)} \times J_C ,\mathbb Q)^{\operatorname{Sp}(V)}$ is also non-degenerate and hence Poincar\'e duality is satisfied. This completes the proof. 
\end{proof}

\begin{theorem}\label{thm:cn-j}
    Let $\xi_{i,g}$ be the pullback of $\xi_g$ over the projection onto the $i$-th factor $\mC_{2,g}^n\times_{\h_g}\jac^g_{2,g}\ri\mC_{2,g}\times_{\h_g}\jac^g_{2,g}$.
For $1\leq n\leq g+6$, the ring $A^*(\mC_{2,g}^n\times_{\h_g}\jac^g_{2,g})$ is multiplicatively generated by $\{\xi_{i,g},\psi_{i},p_\jac^*\T_d,\Delta_{ij}\}_{1\leq i\neq j\leq n}$.
\end{theorem}
\begin{proof}
 From the proof of Theorem \ref{thm:poincare-cnj}, we get \begin{equation*}
      A^* (\mC_{2,g}^{n}\times_{\h_g}\jac^g_{2,g}) \simeq  H^{2*} (C^{n} \times J_C ,\mathbb Q)^{\operatorname{Sp}(V)}.
 \end{equation*}
 The result now follows from observing that the restrictions of $\{\xi_{i,g},\psi_{i},p_\jac^*\T_d,\Delta_{ij}\}_{1\leq i\neq j\leq n}$ give a set of multiplicative generators for $H^{2*}(C^n \times J_C,\mathbb Q)^{\operatorname{Sp}(V)}$.  
\end{proof}
\begin{corollary}
    Let $p_\jac: \mC_{2,g}\times_{\h_g}\J^d_{2,g}\ri \J^d_{2,g}$ be the projection map. We have
    \begin{align*}A^*(\mC_{2,g}\times_{\h_g}\J^d_{2,g})\simeq \Q[a,b ,c ,d]/I\end{align*}
    where $a:= p_\mC^*K, b:= \Xi_d,c:= p_\J^*p_{\J *}(\Xi_d^2),d: = p_\J^*p_{\J *}(\Xi_dp_\mC^*K)$ and the ideal\\
     $I$ is generated by
    $p_\mC^*K^2, \xi_d p_\mC^*K, \xi_d^2+\dfrac{p_\mC^*K.p_\J^*\T_d}{g-1}, p_\J^*\T_d^{g+1},(\xi_d-\dfrac{d}{2g-2}p_\mC^*K)p_\J^*\T_d^{g}$. 
\end{corollary}
\begin{proof}
Applying Theorem \ref{thm:gerbe-ring}, and \cite[Theorem 1.1]{larson-24} or Theorem \ref{simpler proof of Larson}, we have that 
$A^*(\mC_{2,g}\times_{\h_g}\J^d_{2,g})$ is multiplicatively generated by $p_\mC^*K, \xi_d, p_\J^*\ka_{-1,2},p_\J^*\ka_{0,1}$ and the claimed ideal of relations. 
\end{proof}
\begin{theorem}
    For $0\leq n\leq  g+6$ and $0\leq i\leq g$, we have 
    $$\dim A^i(\mC_{2,g}^{(n)}\times_{\h_g}\jac^g_{2,g}) = \#\{(a,b,c)\in \Z_{\geq 0}^3|\; a+b+c\leq i, 2a+b+c\leq n\}.$$
\end{theorem}
\begin{proof}
From the proof of Theorem \ref{thm:poincare-cnj}, it's enough to analyze $H^* (C^{(n)}\times J_C,\mathbb Q)^{\operatorname{Sp(V)}}$. We know $H^* (C^{(n)},\mathbb Q) = \oplus_{r+b \leq n}\bigwedge ^r V$ as graded $\operatorname{Sp}(V)$-representations with $\bigwedge^r V$ in degree $r+2b$. So we get\begin{equation*}
    \begin{aligned}
        \dim_{\mathbb Q} A^i(\mC_{2,g}^{(n)}\times_{\h_g}\jac^g_{2,g}) =\sum_{\substack{r+b\leq n \\ r+s+2b=2i}} \dim_{\mathbb Q} \left (\bigwedge^r V\otimes \bigwedge^{s} V \right )^{\operatorname{Sp}(V)}.
    \end{aligned}
\end{equation*}From \cite[Theorem 5.5.15, Corollary 5.5.16]{MR2522486}, we get \begin{equation*}
    \bigwedge^r V = \bigoplus_{\substack{ c+2a=r \\ 0\leq c\leq g \\ 0\leq a\leq g-c}} P^c
\end{equation*}where $P^c$ is the irreducible primitive representation of highest weight $\omega_c$ and $\dim P^c = \binom{2g}{c} -\binom{2g}{c-2}$. Moreover these $P^c$'s are self dual and pairwise non-isomorphic. So when we write \begin{equation*}
    \begin{aligned}
         \bigwedge^r V  \otimes  \bigwedge^{s} V =   \bigoplus_{\substack{ c+2a=r \\ 0\leq c\leq g \\ 0\leq a\leq g-c}} \bigoplus_{\substack{d+2b=s \\ 0\leq d\leq g \\ 0\leq b\leq g-d}} P^c \otimes P^{d},
    \end{aligned}
\end{equation*}we get \begin{equation*}
    \dim_{\mathbb Q} \left (\bigwedge^r V  \otimes  \bigwedge^{s} V \right ) ^{\operatorname{Sp}(V)} = \# \{c \mid  0\leq c\leq \min\{ r,s,2g-r,2g-s \}, \ c\equiv r\equiv s \mod 2\}.
\end{equation*}
Using the above result we get after a straight-forward but tedious computation that \begin{equation*}
    \begin{aligned}
 \dim_{\mathbb Q} A^i(\mC_{2,g}^{(n)}\times_{\h_g}\jac^g_{2,g})  & = \sum_{\substack{r+b \leq n \\ r+s+2b=2i} } \sum_{\substack{c+2a=r \\ c+2e=s \\ 0 \leq c\leq g\\ 0\leq a,e\leq g-c}} 1 &= \sum_{\substack{a,b,c,e\geq 0 \\ 2a+b+c\leq n \\ a+b+c+e=i \\ c\leq g,  \ a,e \leq g-c}}1 &= \sum_{\substack{a,b,c\geq 0 \\ 2a+b+c\leq n \\ a+b+c\leq i }}1   ,
    \end{aligned}
\end{equation*}where the last equality follows from $i\leq g$, so the sum remains the same. This completes the proof. 
\end{proof}\
We have the following corollaries.
\begin{corollary}\label{cor:c-j}
  The ring $A^*(\mC_{2,g}\times_{\h_g}\jac^d_{2,g})$ is multiplicatively generated by $p_\mC^*K, \xi_d, p_\jac^*\T_d$. In particular, we have
    \begin{align*}
      A^*(\mC_{2,g}\times_{\h_g}\jac^d_{2,g})\simeq \Q[p_\mC^*K, \xi_d, p_\jac^*\T_d] / \langle p_\mC^*K^2, \xi_d p_\mC^*K, \xi_d^2+\dfrac{p_\mC^*K.p_\jac^*\T_d}{g-1}, p_\jac^*\T_d^{g+1},
	  (\xi_d-\dfrac{d}{2g-2}p_\mC^*K)p_\jac^*\T_d^{g} \rangle .
    \end{align*}
\end{corollary}
\begin{corollary}
    Let $p_\jac: \mC_{2,g}\times_{\h_g}\J^d_{2,g}\ri \J^d_{2,g}$ be the projection map. We have
    \begin{align*}A^*(\mC_{2,g}\times_{\h_g}\J^d_{2,g})\simeq \Q[a,b ,c ,d]/I\end{align*}
    where $a:= p_\mC^*K, b:= \Xi_d,c:= p_\J^*p_{\J *}(\Xi_d^2),d: = p_\J^*p_{\J *}(\Xi_dp_\mC^*K)$ and the ideal\\
     $I$ is generated by
    $p_\mC^*K^2, \xi_d p_\mC^*K, \xi_d^2+\dfrac{p_\mC^*K.p_\J^*\T_d}{g-1}, p_\J^*\T_d^{g+1},(\xi_d-\dfrac{d}{2g-2}p_\mC^*K)p_\J^*\T_d^{g}$.
\end{corollary}
\begin{theorem}\label{thm:hilb-stable-denom}
Let $h_U(T)$ be the Hilbert series of $\mcalu^\h(2,d,g)$. We have the following:
\begin{equation} \begin{aligned}
(1-T)^2(1-T^2)^2(1-T^3)h_U(T) =& 1 - T^g - 3T^{g+1} - 3T^{g+2} - T^{g+3}+ (2g+3)T^{2g+1} + (2g+3)T^{2g+2}\\
&- (2g+3)T^{2g+4} - (2g+3)T^{2g+5}+ T^{3g+3} + 3T^{3g+4} + 3T^{3g+5}\\& + T^{3g+6}- T^{4g+6}. \end{aligned}
\end{equation}
\end{theorem}

\begin{proof}

Let $c_{i,n} = \#\{(a,b,c) \in \mathbb{Z}_{\geq 0}^3 \mid a+b+c \leq i, 2a+b+c \leq n\}$. 
We define the formal power series $F(T,y) = \sum_{i,n \geq 0} c_{i,n} T^i y^n$. By summing over the independent variables $a,b,c$ and encoding the cumulative inequalities via $(1-T)^{-1}(1-y)^{-1}$, we obtain the exact rational function:
\begin{equation*}
F(T,y) = \frac{1}{(1-T)(1-y)(1-Ty^2)(1-Ty)^2}.
\end{equation*}
For a fixed $n$, the power series 
$$R_n(T) := \sum_{i=0}^\infty c_{i,n} T^i$$
 is extracted via Cauchy's integral formula over a small contour :
\begin{equation}\label{eq:contour}
R_n(T) = \frac{1}{2\pi i} \oint_{\gamma} F(T,y) \frac{dy}{y^{n+1}} \text{ with } \gamma := \{y \in \mathbb{C} : |y| = \epsilon\}.
\end{equation}

We now construct the palindromic polynomial, by Theorem \ref{thm:poincare-cnj}, $P_{g,n}(T)$ directly from the rational function $R_n(T)$. 
It is easy to see that the sequence of coefficients stabilizes: $c_{i,n} = c_{n,n}$ for all $i \geq n$.

Consequently, the infinite series $R_n(T)$ admits the exact rational representation since $n<g$:
\begin{equation*}
R_n(T) = \sum_{i=0}^{g} c_{i,n} T^i + c_{n,n} \frac{T^{g+1}}{1-T}.
\end{equation*}
The target polynomial $P_{g,n}(T)$ is palindromic of degree $n+g$. 
Because $c_{i,n}$ is constant for $n \leq i \leq g$, the truncation and palindromic reflection can be captured exactly by adding the reflected rational tail. 
Evaluating $T^{n+g} R_n(T^{-1})$, we get:
\begin{equation*}
P_{g,n}(T) = R_n(T) + T^{n+g} R_n(T^{-1}).
\end{equation*}

We now apply this exact rational identity to $h_U(T)$. 
Let the weight function be $W_n(T) := T^n + T^{3g-3-2n}$ for $0\leq n \leq g-2$, and $W_{g-1}(T) := T^{g-1}$. 
We can express $h_U(T)$ entirely in terms of the rational operators:
\begin{equation*}
h_U(T) = \sum_{n=0}^{g-1} W_n(T) P_{g,n}(T) = \sum_{n=0}^{g-1} W_n(T) \left( R_n(T) + T^{n+g} R_n(T^{-1}) \right).
\end{equation*}
By interchanging the sum with the contour integral extraction $R_n(T)$ in (\ref{eq:contour}), we get the following expression for $h_U(T)$:
\begin{align*}
h_U(T) = \frac{1}{2\pi i} \oint_{\gamma} \left[ F(T,y) \sum_{n=0}^{g-1} \frac{W_n(T)}{y^{n+1}} + F(T^{-1}, y) \sum_{n=0}^{g-1} \frac{W_n(T) T^{n+g}}{y^{n+1}} \right] dy.
\end{align*}

By the Global Residue Theorem, evaluating this integral reduces to summing the residues at the exterior poles of $F(T,y)$ and $F(T^{-1},y)$, which are $y=1$, $y=T^{\pm 1}$, and $y=\pm T^{\pm 1/2}$ based on an arbitrary choice of the two squareroots of $T$.
 The integrand is $\omega_1(T,y) + \omega_2(T,y)$ where:
\begin{align*}
\omega_1(T,y) &= \frac{1}{(1-T)(1-y)(1-Ty^2)(1-Ty)^2 y} \sum_{n=0}^{g-1} \frac{W_n(T)}{y^n}, \\
\omega_2(T,y) &= \frac{-T^4}{(1-T)(1-y)(T-y^2)(T-y)^2 y} \sum_{n=0}^{g-1} \frac{W_n(T) T^{n+g}}{y^n}.
\end{align*}

We calculate the explicit residues at the isolated poles of $\omega$:

\textbf{ The simple pole at $y=1$}:
For $y=1$, both $\omega_1$ and $\omega_2$ have simple poles. We get: 
\begin{equation*}
\operatorname{Res}_{y=1} \omega = \frac{-1}{(1-T)^4} \sum_{n=0}^{g-1} W_n(T) -\frac{T^4}{(1-T)^4} \sum_{n=0}^{g-1} W_n(T) T^{n+g}.
\end{equation*}

\textbf{ The double pole at $y=1/T$}:
This is a pole of order 2 for $\omega_1$ only; $\omega_2$ is analytic at $y=1/T$, so its residue is $0$. We get:
\begin{equation*}
\operatorname{Res}_{y=1/T} \omega = \frac{-3T}{(1-T)^4} \sum_{n=0}^{g-1} W_n(T) T^{n+1} - \frac{1}{(1-T)^3} \sum_{n=0}^{g-1} W_n(T) (n+1) T^{n+2}.
\end{equation*}

\textbf{ The double pole at $y=T$}:
This is a pole of order 2 for $\omega_2$ only; $\omega_1$ is analytic at $y=T$, so its residue is $0$. We get:
\begin{equation*}
\operatorname{Res}_{y=T} \omega = \frac{-3T^3}{(1-T)^4} \sum_{n=0}^{g-1} W_n(T) T^{g-1} + \frac{T^3}{(1-T)^3} \sum_{n=0}^{g-1} W_n(T) (n+1) T^{g-2}.
\end{equation*}

\textbf{ The simple poles at $y = \pm 1/\sqrt{T}$}:
These are simple poles of $\omega_1$ from the factor $(1-Ty^2)$. The roots are $y_0 = \pm 1/\sqrt{T}$. We get:
\begin{align*}
\operatorname{Res}_{y=1/\sqrt{T}} \omega &= \frac{-1}{2(1-T)\left(1-\frac{1}{\sqrt{T}}\right)(1-\sqrt{T})^2} \sum_{n=0}^{g-1} W_n(T) T^{n/2}, \\
\operatorname{Res}_{y=-1/\sqrt{T}} \omega &= \frac{-1}{2(1-T)\left(1+\frac{1}{\sqrt{T}}\right)(1+\sqrt{T})^2} \sum_{n=0}^{g-1} W_n(T) (-1)^n T^{n/2}.
\end{align*}

\textbf{ The simple poles at $y = \pm \sqrt{T}$}:
These are simple poles of $\omega_2$ from the factor $(T-y^2)$. The roots are $y_0 = \pm \sqrt{T}$. We get:
\begin{align*}
\operatorname{Res}_{y=\sqrt{T}} \omega &= \frac{T^3}{2(1-T)(1-\sqrt{T})(T-\sqrt{T})^2} \sum_{n=0}^{g-1} W_n(T) T^{n+g} T^{-n/2}, \\
\operatorname{Res}_{y=-\sqrt{T}} \omega &= \frac{T^3}{2(1-T)(1+\sqrt{T})(T+\sqrt{T})^2} \sum_{n=0}^{g-1} W_n(T) (-1)^n T^{n+g} T^{-n/2}.
\end{align*}
Set
\begin{equation*} \begin{aligned}
&X :=\sum_{n=0}^{g-1} W_n =  \frac{1-T^g}{1-T}+\frac{T^{g+1}-T^{3g-1}}{1-T^2},\end{aligned}\end{equation*} \begin{equation*} \begin{aligned}
&Y :=\sum_{n=0}^{g-1}W_nT^n= \frac{1-T^{2g}}{1-T^2}+\frac{T^{2g-1}-T^{3g-2}}{1-T}, \end{aligned} \end{equation*} \begin{equation*} \begin{aligned}
Z: = &\sum_{n=0}^{g-1}W_n T^{n/2} = \frac{1-T^{3g-3/2}}{1-T^{3/2}}, \end{aligned}\end{equation*} \begin{equation*} \begin{aligned}
    Z' := &\sum_{n=0}^{g-1} W_n(-1)^nT^{n/2} = \frac{1+T^{3g-3/2}}{1+T^{3/2}},
    \end{aligned}
\end{equation*}
\begin{equation*} \begin{aligned}
    V := &\sum_{n=0}^{g-1}(n+1)W_nT^n = \frac{gT^{2g+2}+1-(g+1)T^{2g}}{(1-T^2)^2}+T^{2g-1}\frac{(g-1)+T^g-gT}{(1-T)^2},
    \end{aligned}
\end{equation*}
\begin{equation*}
\begin{aligned}
    V' := &\sum_{n=0}^{g-1} (n+1)W_n =\frac{gT^{g+1}+1-(g+1)T^g}{(1-T)^2}+\frac{(g-1)T^{g+1}+T^{3g+1}-gT^{g+3}}{(1-T^2)^2}.
    \end{aligned}
\end{equation*}
\noindent Adding the seven listed residue values gives us $-h_U(T)$. More precisely, we have  
\begin{equation*}\begin{aligned}h_U(T)= & X\frac{1+3T^{g+2}}{(1-T)^4}+Y\frac{T^{g+4}+3T^2}{(1-T)^4}+V\frac{T^2}{(1-T)^3}\\ &-V'\frac{T^{g+1}}{(1-T)^3}-Z\frac{T^{g+2}+\sqrt{T}}{2(1-T)(1-\sqrt{T})^3}+Z'\frac{\sqrt{T}-T^{g+2}}{2(1-T)(1+\sqrt{T})^3}.\end{aligned}\end{equation*}
Consequently, we have the claimed expression for $(1-T)^2(1-T^2)^2(1-T^3)h_U(T)$. 
The associated computations are executed via SymPy and the code is attached in Figure \ref{hilb:code}.
\end{proof}
\subsection{Final computations and a numerical proof of Conjecture \ref{conj:hyp}}
We give a numerical proof of surjectivity and hence of Conjecture \ref{conj:hyp}.
\begin{proof}[Alternative proof of Conjecture \ref{conj:hyp}]
By Theorems \ref{thm:hilb-unstable} and \ref{thm:hilb-stable-denom}, the sum of codimension-twisted Hilbert series over HN-strata is given by:
\begin{equation*} \begin{aligned}
    &\dfrac{\hilb (\mcalu^\h(2,d,g))}{1-T}+T^g\hilb \left(\bigsqcup_{n\geq \lfloor d/2\rfloor +1}\ext^{(n,d)}_{(1,2)}(g)\right)  \\ =
&\dfrac{ 1 - T^g - 3T^{g+1} - 3T^{g+2} - T^{g+3}+ (2g+3)\left (T^{2g+1} +T^{2g+2}- T^{2g+4} - T^{2g+5} \right )+ T^{3g+3}}{(1-T)^3(1-T^2)^2(1-T^3)}\\
&+\frac{  3T^{3g+4} + 3T^{3g+5} + T^{3g+6}- T^{4g+6}}{(1-T)^3(1-T^2)^2(1-T^3)}+T^g\dfrac{{g+2\choose 2}T^g+\sum_{1\leq k\leq g}{g-k+2\choose 2}(T^{g-k}+T^{g+k})}{(1-T)^2(1-T^2)}. \end{aligned}
\end{equation*}
Multiplying by the denominator in Theorem \ref{thm:hilb-ka}, $(1-T)^3(1-T^2)^2(1-T^3)$, we get:
$$-T^{g+1}- T^{g+2} +T^{3g+4} +T^{ 3g+5} -T^{ 4g+6} +1 = (1-T^{g+1})(1-T^{g+2})(1-T^{2g+3})$$
which matches the numerator in Theorem \ref{thm:hilb-ka}. The computations are attached in Figure \ref{conj:code}. That is, we have
\begin{equation*}\begin{aligned}(1-T)^3(1-T^2)^2(1-T^3)&\bigg(\dfrac{\hilb (\mcalu^\h(2,d,g))}{1-T}+ T^g\hilb \left(\sqcup_{n\geq \lfloor d/2\rfloor +1}\ext^{(n,d)}_{(1,2)}(g)\right)\\ &-\hilb(H^*_\ka(Bun_{2,d}C) )\bigg)= 0.\end{aligned}\end{equation*}

\noindent The conjecture follows from a straightforward adaptation of \cite[Theorem 4.1]{shubham-r2-ss}.
\end{proof}

\section{Even degree bundles}\label{sec:even-deg}
It remains to prove the conjecture for $\ve^\h(2,d,g)$ for $d$ even. 
Throughout this section, we shall use the fact that the pullback of Chow rings along finite surjective maps of smooth DM-stacks is injective without explicit mention.
\begin{lemma} \label{stratification generation}
    Let $X$ be a smooth DM stack and $R\subset A^* (X)$ be a sub-ring. Suppose $F_0\subset X$ is a smooth locally closed substack admitting a filtration $F_0 \supset F_1\supset\cdots \supset F_r$ by smooth closed substacks of $F_0$. Set $W_i =F_i \backslash F_{i+1}$ for $0\leq i< r$ and $W_r=F_r$. Assume the following:\begin{enumerate}
        \item For every $0\leq i\leq r$, the restriction map $R\rightarrow A^* (W_i)$ is surjective.
        \item $[F_{i+1}]\in A^*(F_i)$ is the restriction of some $a_i \in R$. 
    \end{enumerate}
    Then the restriction map $R\rightarrow A^*(F_0)$ is surjective. 
\end{lemma}
\begin{proof}
    We proceed by descending induction on $i$. For $i=r$, this is (1) since $W_r=F_r$. Assume $R\rightarrow A^* (F_{i+1})$ is surjective. By the excision exact sequence $$ A^* (F_{i+1})\rightarrow A^* (F_i) \rightarrow A^* (W_i) \rightarrow 0,$$ and (1), given $z\in A^* (F_i)$ we can find $\tilde z\in R$ such that $z-\tilde z|_{F_i} = j_* z'$ for some $z'\in A^* (F_{i+1})$ where $j: F_{i+1} \hookrightarrow F_i$ is the inclusion. But by induction hypothesis, $z' =\tilde z'|_{F_{i+1}}$ for some $\tilde z' \in R$ so we have $j_* (z') =j_* (\tilde z'|_{F_{i+1}}) = (\tilde z'|_{F_i}). j_*([F_{i+1}]) = (\tilde z'. a_i)|_{F_i}$. This completes the inductive step and hence the proof. 
\end{proof}
\begin{lemma}\label{H surjects onto Z}
      Let \begin{equation*}
        \mathcal Z_{g,m}:= \left \{ (C,p_1,\dots, p_m)\in \mathcal H_{g,m}\mid p_1\text{ is a Weierstrass point.} \right \}.
    \end{equation*}
We have the closed immersion $j_m: \mathcal Z_{g,m} \hookrightarrow \mathcal H_{g,m}$. Then $$j_m^* :A ^* (\mathcal H_{g,m}) \rightarrow A^* (\mathcal Z_{g,m})$$ is surjective for $m\leq 2g+5$.
\end{lemma}
\begin{proof}
    We adopt the strategy in \cite{cl-hyp}. Let $\mathcal H^o_{g,m}$ (respectively $\mathcal Z^o_{g,m}$) denote the open locus where no two marked points are conjugate under the hyper-elliptic involution. Set $W_i^{sm}:= \big \{(C, p_1,\dots, p_m) \in \mathcal Z^o_{g,m} \mid p_1,\dots, p_i \text{ are Weierstrass, }p_{i+1}\text{ is not Weierstrass.}   \big \}$ \footnote{Our $W_i^{sm}$'s are open in the $W_i$'s of \cite{cl-hyp} because we want $C$ to be smooth.}. In the notation of \cite{cl-hyp}, we have the Cartesian diagram \begin{tikzcd}
    W^{sm}_i\arrow{d} \arrow{r}{\subset_{\text{open}}} &  W_i \arrow{d} \\ \mathcal H^o_{g,m} \arrow{r}{\subset_{\text{open}}}& \mathcal I^o_{g,m} 
\end{tikzcd}.  From \cite[Lemma 4.6]{cl-hyp}, we know $A^* (\mathcal I^0_{g,m}) \rightarrow A^* (W_i)$ is surjective for $1\leq i<m \leq g$ and hence by excision, we know $A^* (\mathcal H_{g,m}) \rightarrow A^* (W^{sm}_i)$ is surjective for $1\leq i<m\leq g$. From \cite[Lemma 3.8]{cl-hyp}, we know that $A^* (W^{sm}_m) =\mathbb Q$ and hence $A^*(\mathcal H_{g,m}) \ri A^*(W^{sm}_m)$ is trivially surjective. Moreover, if we set $$F^{sm}_i:= \left \{ (C,p_1,\dots, p_m) \in \mathcal Z^o_{g,m} \mid p_1,\dots, p_i\text{ are Weierstrass points.} \right \}, $$
then $F^{sm}_i \backslash F^{sm}_{i+1} = W^{sm}_i$. From \cite[Lemma 2.9]{cl-hyp}, we see that $[F^{sm}_i]$ is the restriction of a product of divisors. So we can apply \ref{stratification generation} to conclude that the restriction map $A^* (\mathcal H_{g,m}) \rightarrow A^* (\mathcal Z^o_{g,m})$ is surjective for $1\leq m\leq g$.  We now prove $A^* (\mathcal H_{g,m}) \rightarrow A^* (\mathcal Z_{g,m})$ is surjective for $1\leq m \leq g$. For $m=1$, $\mathcal Z_{g,1} =\mathcal H^w_g$ and we know from \cite[Lemma 4.1]{clader2025chowringshurwitzspaces} that $A^* (\mathcal H_g^w) =\mathbb Q$ so $A^* (\mathcal H_{g,1}) \rightarrow A^* (\mathcal Z_{g,1})$ is trivially surjective\footnote{In the notation of \cite{clader2025chowringshurwitzspaces}, we have $\mathcal H^w_g =\mathcal H_{2,g}(2)$.}. We also note that this implies $\mathcal Z_{g,1}$ is irreducible and considering the forgetful map $\mathcal Z_{g,m} \rightarrow \mathcal Z_{g,1}$, we can deduce that each $Z_{g,m}$ is irreducible. We proceed by induction on $m$. Suppose we have proved $A^* (\mathcal H_{g,m-1}) \rightarrow A^* (\mathcal Z_{g,m-1})$ is surjective. Set $E_{i,j}:=\mathcal Z_{g,m}\cap D_{i,j}$ for $2\leq i< j\leq m$. Forgetting $p_j$ identifies $E_{i,j}$ with an open sub-stack of $\mathcal Z_{g,m-1}$ and in particular, we note that $E_{i,j}$ is also irreducible. Thus from the commutative diagram \begin{tikzcd}
    \mathcal H_{g,m}\arrow{r}{\operatorname{Forget}_j} & \mathcal H_{g,m-1}  \\ E_{ij}\arrow{u} \arrow{r} & V_{i,j}\arrow{r}{ \subset_{\text{open}}}\arrow{u}& \mathcal Z_{g,m-1} \arrow{ul}
\end{tikzcd}, and the induction hypothesis, we get the restriction map $A^* (\mathcal H_{g,m}) \rightarrow A^* (E_{i,j})$ is surjective for $2\leq i< j\leq m$. Now since $E_{i,j}$ is irreducible, restriction of $D_{i,j}$ $\in \mathbb Q_{>0}. [E_{i,j}]$ and in particular, we have $[E_{i,j}]\in A^*(\mathcal Z_{g,m})$ is the restriction of a class in $A^* (\mathcal H_{g,m})$. A similar argument as in \ref{stratification generation} then shows us the restriction $A^* (\mathcal H_{g,m}) \rightarrow A_* (\mathcal Z_{g,m})$ is surjective for $1\leq m \leq g$ in view of the excision exact sequence $$ \bigoplus_{2\leq i< j\leq m} A_* (E_{i,j}) \rightarrow A_* (\mathcal Z_{g,m}) \rightarrow A^* (\mathcal Z^o_{g,m})\rightarrow 0.$$ 

\noindent We now tackle the case $g+1\leq m \leq 2g+5$. We use the well-known fact that a smooth hyper-elliptic curve $C$ admits an embedding  into $\mathbb P^1\times \mathbb P^1$ as $C= \left \{ F(x,y) =0 \right \}$, where $F$ has bi-degree $(g+1,2)$ and the hyper-elliptic map is the projection onto the first factor. From \cite{cl-hyp}, we know $\mathcal H^o_{g,m}$ is an open substack of $\mathbb P(\mathcal F)$ where $\mathcal F$ is a bundle of rank $3g+7-m$ (see \cite[Equation 4.11]{cl-hyp}) on $\mathcal B^o_{g,m}$ where $\mathcal B^o_{g,m}$ is the open locus in the configuration stack $\mathcal B_{g,m}$ consisting of points $(q_1,\dots, q_g, q_{g+2},\dots,q_{m})\in \mathbb P^1\times \mathbb P^1$ such that the $m-1$ points impose linearly independent conditions on equations of bi-degree $(g+1,2)$. Let $\mathcal B^{o,w}_{g,m}$ denote the open locus consisting of points $(q_1=(x_1,y_1),\dots ,q_g=(x_g,y_g), q_{g+2}=(x_{g+2},y_{g+2}),\dots q_m =(x_m,y_m))$ such that the $m-1$ points and the condition $\partial_y F(q_1)=0$ impose  linearly independent conditions on equations of bi-degree $(g+1,2)$. Clearly, $\mathcal B^{o,w}_{g,m} \subset \mathcal B^{o}_{g,m}$. For $b\in \mathcal B^{o,w}_{g,m}$ as above consider the one dimensional vector space $L_b=\operatorname{Hom} \left (T_{y_1}\mathbb P^1, \mathcal O(g+1,2)|_{q_1} \right )$. As $b$ varies, this gives us a line bundle $\mathcal L$ on $\mathcal B^{o,w}_{g,m}$.  Moreover we have $\mathcal F \xrightarrow{d_y} \mathcal L$, which is fiber-wise given by $F\mapsto \partial_y F(q_1)$. By definition of $\mathcal B^{o,w}_{g,m}$ the kernel of this morphism is a vector bundle $\mathcal F^w$ on $\mathcal B^{o,w}_{g,m}$ of rank $\operatorname{rk}\mathcal F-1=3g+6-m$. We claim that $\mathcal Z^o_{g,m} \rightarrow \mathbb P(\mathcal F) \rightarrow \mathcal B^o_{g,m}$ factors through $\mathcal B^{o,w}_{g,m}$. To see this, consider the divisor $\displaystyle{\Gamma:=2p_1 +\sum_{\substack{j=2 \\ j\neq g+1}}^m p_j}$. By \cite[Lemma 4.2]{cl-hyp}, the evaluation map \begin{equation*}
    H^0(\mathbb P^1\times \mathbb P^1 , \mathcal O(g+1,2)) \rightarrow H^0( \Gamma, \mathcal O_{\Gamma}(g+1,2) ) 
\end{equation*}is surjective since $m\leq 2g+5$. In other words, the $m-1$ points and the vanishing of the derivative impose linearly independent conditions on equations of bi-degree $(g+1,2)$. So we get the image of $\mathcal Z^o_{g,m}$ lands in $\mathcal B^{o,w}_{g,m}$. 

Now consider the Cartesian square \begin{tikzcd}
    \mathcal H^o_{g,m}\times_{\mathbb P(\mathcal F)} \mathbb P(\mathcal F^w)\arrow{d} \arrow{r}[swap]{\subset_{\text{open}}} & \mathbb P(\mathcal F^w) \arrow{d} \\ \mathcal H^o_{g,m} \arrow{r}{\subset_{\text{open}}}& \mathbb P(\mathcal F)
\end{tikzcd}. We know the top horizontal arrow is an open immersion and it parametrizes precisely the points $\{ (C,p_1,p_2,\dots p_m) \in \mathcal H^o_{g,m}\mid p_1\text { is Weierstrass.}\}$. So we have $  \mathcal H^o_{g,m}\times_{\mathbb P(\mathcal F)} \mathbb P(\mathcal F^w) =\mathcal Z^o_{g,m}$. Now the restriction map $A^* (\mathbb P(\mathcal F)) \rightarrow A^* (\mathbb P(\mathcal F^w))$ is surjective and hence so is the restriction map $A^* (\mathcal H_{g,m})\rightarrow A^*(\mathcal H^o_{g,m}) \rightarrow A^* (\mathcal Z^o_{g,m})$. To show that $A^* (\mathcal H_{g,m})\rightarrow A^* (\mathcal Z_{g,m})$ is surjective, we use the same argument as in $1\leq m \leq g$ case.

\end{proof}
\begin{corollary}\label{chow of Z}
    We have $A^*(\mathcal Z_{g,m})=\mathbb Q$ for $1\leq m\leq 2g+5$. 
\end{corollary}
\begin{proof}
    We know from \cite{cl-hyp}, that for $1\leq m\leq 2g+5$, $A^* (\mathcal H_{g,m}) =\mathbb Q[\psi_1,\dots, \psi_m]/ (\psi_1,\dots,\psi_m)^2$. Let $\mathcal C^{\mathcal Z}_{g,m} \xrightarrow{\pi} \mathcal Z_{g,m}$ denote the universal curve and $s_1,\dots, s_m$ are the sections. Since $s_1$ is a Weierstrass section, we have $\omega_\pi= \mathcal O((2g-2)s_1) \otimes \pi^* M$ for some line bundle $M$ on $\mathcal Z_{g,1} = \mathcal H^w_{g}$. For $i>1$, $s_1, s_i$ are disjoint hence $\psi_i =c_1(s_i^*\omega_\pi)=s_i^*\pi^* c_1( M) =0$ since $A^1(\mathcal H^w_g)=0.$ For $i=1$, we have $\psi_1 $ is pulled back from a divisor in $\mathcal H^w_g$ and hence $0$ since $A^1(\mathcal H^w_g)=0$. This completes the proof. 
\end{proof}
\begin{theorem}\label{thm:weier-pull}
    For all $g\geq 2, 0\leq n\leq 2g+4$, the pullback map $A^*(\h_{g,n})\ri A^*(\h_{g,n}^w)$ is an isomorphism.
\end{theorem}
\begin{proof}
    Let $q_n :\mathcal H^w_{g,n} \rightarrow \mathcal H_{g,n}$ denote the finite etale map of degree $2g+2$. It is clear that $q^*_n : A^* (\mathcal H^w_{g,n}) \rightarrow A^*(\mathcal H_{g,n})$ is injective. So we only need to establish surjectivity. 

\noindent For $1\leq i \leq n$, set $D_i = \{w=p_i \}\subset \mathcal H^w_{g,n}$. Note that $D_i \simeq \mathcal Z_{g,n}$ and hence $A^* (D_i)=\mathbb Q$. Now if we set $U: = \mathcal H^w_{g,n}\backslash \bigcup_{i=1}^n D_i$, then $U\simeq \mathcal Z_{g,n+1}$. So we have for $1\leq n\leq 2g+4$, $A^*(U) \simeq \mathbb Q$ by Corollary \ref{chow of Z}. It follows from a straightforward argument involving the excision sequence that \begin{equation*}
    A^* (\mathcal H^w_{g,n}) =\operatorname{Span}_{\mathbb Q} \left \{ 1, [D_1],\dots [D_n] \right \}. 
\end{equation*}
But if $\mathcal C^w_{n,g} \xrightarrow{\pi} \mathcal H^w_{n,g} $ is the universal curve, we have $\omega_\pi = \mathcal O((2g-2)w)$ upto tensoring by a line bundle pulled back from $\mathcal H^w_g$. It follows that  $q_n^* \psi_i= c_1(s_i^* \omega_\pi) =(2g-2)s_i^* c_1(\mathcal O(w))= (2g-2)c_1(\mathcal O_{\mathcal H^{w}_{g,n}}(D_i))=(2g-2)[D_i]$. This shows $q_n^*$ is surjective.

\end{proof}

\begin{corollary}\label{cor:hyp-weier-powers}
    The pullback map $A^*(\mC_{2,g}^n)\ri A^*(\times_{\h_{g}^w}^n\mC_{2,g}^w)$ is an isomorphism for all $g\geq 2, 0\leq n\leq 2g+4$.
\end{corollary}
\begin{proof}
Using the stratification in Lemma \ref{lem:gens} and applying Theorem \ref{thm:weier-pull}, we have surjectivity of the pullback map.
\end{proof}
\begin{lemma}\label{lem:pull-iso}
    Let $\mathcal{X},\x',\y,\y'$ be smooth DM-stacks of finite type such that we have the following Cartesian diagram:
    $$\begin{tikzcd}
\x' \arrow[r, "f'"] \arrow[d, "g'"]\ar[dr, phantom, "\square"] & \y' \arrow[d, "g"] \\
\x \arrow[r, "f"] & \y
    \end{tikzcd}$$
    Additionally, assume that $f,f'$ are proper and generically finite, $g,g'$ are finite and surjective, and that the pullback
    $g'^*:A^*(\x)\ri A^*(\x')$ is an isomorphism. Then the pullback $g^*:A^*(\y)\ri A^*(\y')$ is also an isomorphism.
\end{lemma}
\begin{proof}
    We use the formula $g^*f_* = g'^{*}f'_*$ as morphisms $A_*(\x)\ri A^*(\y')$. Surjectivity of $g^*$ follows from the surjectivity of $g'^*, f'_*$.
\end{proof}
\begin{definition}
    Let $\ve_w^\h(2,d,g)$ be the universal moduli stack of rank $2$, degree $d$ vector bundles on hyperelliptic curves of genus $g$ with a marked Weierstrass point.
\end{definition}
\begin{theorem}\label{thm:vec-pull}
    The pullback map $A^*(\ve^\h(2,d,g))\ri A^*(\ve_w^\h(2,d,g))$ is an isomorphism for all $g\geq 2$ and $d$ odd.
\end{theorem}
\begin{proof}
  Following Theorem \ref{thm:sfp} and applying Corollary \ref{cor:hyp-weier-powers} and Lemma \ref{lem:pull-iso} that the pullbacks
$$A^*(\mC_{2,g}^{(n)}\times_{\h_g} \jac^{d'}_{2,g})\ri A^*(\mC_{2,g}^{(n)}\times_{\h_g} \jac^{d'}_{2,g}\times_{\h_g}\h_g^w), 
A^*(\jac^{d'}_{2,g}\times_{\h_g}\jac^{d"}_{2,g})\ri A^*(\jac^{d'}_{2,g}\times_{\h_g}\jac^{d"}_{2,g}\times_{\h_g}\h_g^w)$$
are isomorphisms for all $0\leq n\leq g+4$ and $d',d''\in\Z$.

The claim now follows from Corollary \ref{cor:decomp-ve} and Theorem \ref{thm:sfp-u}.
\end{proof}
\begin{theorem}\label{thm:conj-even}
    Conjecture \ref{conj:hyp} holds for all $g\geq 2$ and $d$ even.
\end{theorem}
\begin{proof}
For $e\in\Z$, let $\p_{e,g}:= \p_{\ve_w^\h(2,e,g)}(\e^w_{e,g}|_{w\times_{\h_g^w}\ve_w^\h(2,e,g)})$ where $w:\h_g^w\ri\mC_{2,g}^w$ is the Weierstrass section and $\e^w_{e,g}\ri \mC_{2,g}^w\times_{\h_g^w}\ve^\h_w(2,e,g)$ is the Poincar\'e bundle.

By \cite[Theorem 4.2]{shubham-r2-ss}, we have that $\p_{e,g}\ri \ve_w^\h(2,e,g)$ is a $\p^1$-bundle and $\p_{e,g}\simeq \p_{e-1,g}$ for all $e\in\Z$.
Applying the projective bundle formula, we get that 
$$\hilb (\ve_w^\h(2,d,g))(1+T) = \hilb(\p_{d,g}) = \hilb(\p_{d-1,g}) = \hilb(\ve^\h_w(2,d-1,g))(1+T)$$ for all $g\geq 2$ and $d$ even.

\noindent Applying Theorem \ref{thm:vec-pull}, we get that $\hilb (\ve_w^\h(2,d,g)) = \hilb(\ve^\h(2,d-1,g))$ for all $g\geq 2$ and $d$ even.
The claim follows from Theorem \ref{thm:conj-odd} and this also gives a proof of Theorem \ref{global Chow = monodromy invariant of a fiber} in the even case.
\end{proof}
\begin{rem}
   Note that the proof of Theorem \ref{thm:conj-even} does not depend on the explicit computations of the Hilbert series for the odd case.
\end{rem}
\begin{rem}
    We also note here that our proof strategy combined with the motivic decomposition of the moduli space of stable bundles of rank $3$ and degree $d$ with $\gcd (d,3)=1$ over a smooth curve as in \cite{MR4399901}, gives us analogues of Conjecture \ref{conj:hyp} for $\mathcal U^\mathcal H(3,d,g)$ and $\mathscr U^\mathcal H(3,d,g)$ for $g\leq 8$, and $\mathcal Vec^\mathcal H(3,d,g)$ for $g\leq 6$. 
\end{rem}

\section{Remarks on the higher rank cases}\label{higher rank section}
In this section, we indicate how our strategy can be used to prove some partial results towards Conjecture \ref{conj:hyp} or Theorem \ref{global Chow = monodromy invariant of a fiber} for higher rank stable bundles. Let $r,d$ be co-prime integers with $r\geq 1$. We start with del Bano's formula in $\widehat K_0 ({\mathrm {M}}(k))$ for an arbitrary field $k$ of $\operatorname{char}0$. From \cite[Theorem 4.11]{del-bano}, we can write \begin{equation*}
    [U(r,d.g)] = \dfrac{F_{r,d,g}\left(\mathbb L, [Jac(C)],[Sym^1C],\dots, [Sym^{g-1}C] \right)}{\prod_{m=1}^r(1-\mathbb L^m)^{N_m}} \in \widehat K_0 (\mathrm M(k)) 
\end{equation*}  where $F_{r,d,g}$ is a polynomial with integer coefficients. Further from the expression in loc. cit., if we assign $\mathbb L$ has $\deg 0$, $Jac(C)$ has degree $g$ and $Sym^i C$ has degree $i$, then we can choose an expression where $F_{r,d,g}$ has degree $\leq rg$. Clearing denominators and re-arranging terms with negative coefficients, we get 
\begin{equation}
    \begin{aligned}
    &[U(r,d,g)] (1+R_{r,d,g}(\mathbb L)) +Q_{r,d,g} (\mathbb L, [Jac(C)],[Sym^1C],\dots, [Sym^{g-1}C]) 
    \\ = & P_{r,d,g}(\mathbb L, [Jac(C)],[Sym^1C],\dots, [Sym^{g-1}C]) +[U(r,d,g)] (E_{r,d,g}(\mathbb L)) \ \text{in }\widehat K_0 (\mathrm M(k)) \text{ with }E_{r,d,g}(0) =0 \end{aligned} 
\end{equation} where $P,Q,R$ are polynomials with non-negative integer coefficients. Further each $P,Q$ has $\deg \leq rg $.  

But each side is in the image of the functor $\mathrm M(k) \rightarrow \widehat K_0(\mathrm M(k))$, so thanks to \cite[Proposition 2.2]{MR4399901}, we get an isomorphism in $\mathrm M(k)$ of the form \begin{equation*}
    \mathfrak h(U(r,d,g)) \oplus X_{r,d,g} \simeq Y_{r,d,g} \oplus \mathfrak h(U(r,d,g)) \otimes \mathfrak E_{r,d,g}
\end{equation*}   where $Y_{r,d,g}$ is a direct sum of motives built out of $Jac(C)^k \times \prod_{j=1}^b (Sym^{i_j}C)^{k_j}$ with $kg +\sum_{j=1}^bi_jk_j\leq rg$ and $\mathfrak E_{r,d,g}$ is built out of direct sums of powers of $\mathbb L$ with strictly positive Tate twists.
So we get \begin{equation*}
    \mathfrak h(U(r,d,g)) \hookrightarrow  Y_{r,d,g} \oplus \mathfrak h(U(r,d,g)) \otimes \mathfrak E_{r,d,g} \text{ is a direct summand.}
\end{equation*}Now iterating this, we get \begin{equation*}
    \mathfrak h(U(r,d,g)) \hookrightarrow  \left (\bigoplus_{j=0}^{N-1} Y_{r,d,g} \otimes \mathfrak E_{r,d,g}^{\otimes j} \right ) \bigoplus  \mathfrak h(U(r,d,g)) \otimes \mathfrak E_{r,d,g}^{\otimes N} \text{ is a direct summand}. 
\end{equation*} But $\operatorname{Hom}_{\mathrm M(k)} \left (   \mathfrak h(U(r,d,g)), \mathfrak h(U(r,d,g)) \otimes \mathfrak E_{r,d,g}^{\otimes N}\right ) =0$ for $N>>0$ because of dimension reasons. \begin{equation*}
    \implies  \mathfrak h(U(r,d,g)) \hookrightarrow  \left (\bigoplus_{j=0}^{N-1} Y_{r,d,g} \otimes \mathfrak E_{r,d,g}^{\otimes j} \right ) \text{ is a direct summand}. 
\end{equation*}We indicate how the same thing holds true for smooth curve fibrations. 

\noindent Let $Y^B_{r,d,g}$ denote the relative versions of the motives we get for a single smooth curve.   \begin{equation*} \begin{tikzcd}
    \mathfrak h_B ( U_{\mathcal C/B}(r,d,g)) \arrow[hook, bend left]{r}{\Gamma} &  \bigoplus_{j=0}^{N-1} Y_{r,d,g}^B \otimes \mathfrak E_{r,d,g}^{\otimes j} \arrow[bend left]{l}{\Psi} \end{tikzcd}
\end{equation*} Using the result at the generic point and spreading out we get relative correspondences $\Gamma$, $\Psi$ as above such that $\Psi_\eta \circ \Gamma_\eta =Id$. So we have $Id-\Psi\circ \Gamma$ is nilpotent thanks to Guletskii's Theorem \ref{relative-nilpotence-theorem} since everything is in $\mathrm A(B)$. So $\Gamma$ is a split injection or in other words the relative motive is a direct summand as well.

\noindent By the Edidin-Graham approximation argument used in Theorem \ref{thm:sfp-u}, we get that \begin{equation*}
    \mathfrak h_{\mathcal H_g} (\mathcal U^\mathcal H(r,d,g)) \hookrightarrow \bigoplus_{\alpha} \mathfrak h_{\mathcal H_g}  \left ( (\mathcal J_{2,g}^d)^{k_\alpha} \times_{\mathcal H_g}  (\mathcal C_{2,g}^{(i_{1,\alpha}) })^{\times_{\mathcal H_g}k_{1,\alpha}} \times _{\mathcal H_g} \dots  \times_{\mathcal H_g}(\mathcal C_{2,g}^{(i_{b,\alpha}) })^{\times_{\mathcal H_g}k_{b,\alpha}}   \right )(t_\alpha)
\end{equation*} is a direct summand with $k_\alpha g_\alpha + \sum_{j=1}^{b_{\alpha }} i_{j,\alpha}k_{j,\alpha} \leq rg $. If $rg\leq 2g+6$, it follows from Theorem \ref{thm:c-powers} and the argument of Theorem \ref{thm:sfp} that $\mathcal U^\mathcal H(r,d,g) \rightarrow \mathcal H_g$ satisfies SFP whenever $rg\leq 2g+6$. This completes the proof of Theorem \ref{higher rank}.

\noindent One can then use our monodromy argument to get explicit isomorphisms with the monodromy invariant cohomology ring for such moduli stacks.

\appendix
\vspace{-0.1in}
\section{Applications to multiplicative Chow-Kunneth decompositions}\label{MCK Appendix}
In this section, we indicate how our Strong Franchetta Property can be used to establish existence of multiplicative Chow-Kunneth decompositions. The existence of multiplicative Chow-Kunneth decompositions has been studied extensively in \cite{MR4285110}, \cite{MR4537835}, \cite{MR4803347} etc.  We start with the definition of a Chow-Kunneth decomposition as introduced in \cite{MR1225267}, \cite{MR1225268}.

\begin{definition}
    Let $X$ be a smooth projective variety of dimension $d$. We say that $X$ admits a Chow-Kunneth (CK) decomposition if there exists a decomposition of the diagonal \begin{equation*}
        \Delta_X =\sum_{i=0}^{2d} \pi^i_X \text{ in }A^d(X\times X)
    \end{equation*}such that $\pi^i_X$'s are mutually orthogonal idempotents and $(\pi^i_X)_* H^*(X,\mathbb Q) =H^{i}(X,\mathbb Q)$. A multiplicative Chow-Kunneth(MCK) decomposition is defined as a CK decomposition $\{ \pi^i_X \}_{i=0}^{2d}$ of $X$ that is multiplicative, i.e. it satisfies \begin{equation*}
        \pi^k_X \circ \Delta^{sm} \circ  \left ( \pi^i_X \times \pi^j_X \right ) =0 \text{ in }A^{2d}(X\times \times X)\text{ for }k\neq i+j
    \end{equation*}where $\Delta^{sm}$ is the small diagonal in $X\times X\times X$. 
\end{definition}

Analogously one can define a relative multiplicative Chow-Kunneth Decomposition for a smooth proper family $\mathcal X/B$. 

Let $ U_C (2,\mathcal L)$ denote the moduli space of rank $2$ bundles with $\det =\mathcal L$ with $\deg \mathcal L$ over a hyper-elliptic curve $C$. Note that $U_C(2, \mathcal L)$ is a smooth projective variety of dimension $3g-3$. We shall be primarily interested when $U_C (2,\mathcal L)$ admits a multiplicative Chow-Kunneth decomposition.

\noindent We start with the following easy observation.
\begin{lemma}\label{triple SFP-MCK Lemma}
    Let $U/B$ be a smooth proper family with a relative motivic decomposition in $\mathrm M(B)$ 
    \begin{equation*}
        \mathfrak h_B (U) \simeq \bigoplus_{\alpha } \mathfrak h_B (T_\alpha )(t_\alpha)
    \end{equation*}where $T_\alpha /B$ are smooth proper families and $t_\alpha$'s are the Tate twists. Further assume that \begin{itemize}
        \item $T_\alpha/ B$ admit relative CK decompositions for every $\alpha$, and
        \item For every $\alpha_1, \alpha_2,\alpha_3$, the family $T_{\alpha_1} \times_B T_{\alpha_2} \times_B T_{\alpha_3}\rightarrow B$ satisfies SFP.
    \end{itemize}  Then $U/B$ admits a relative MCK and in particular for every $b\in B$, a geometric point the smooth proper variety $U_b$ admits a MCK decomposition.
\end{lemma}
\begin{proof}
    By functoriality of a motivic isomorphism with the cycle class map, a relative CK on $T_\alpha \rightarrow B$ for every $\alpha$ induces a relative CK decomposition $\{ \pi^i_U \}_{i=0}^{2\dim U-2\dim B }$ for the family $U\rightarrow B$. By the second condition, we know that $U\times_B U \times_B U \rightarrow B$ satisfies SFP. Since we know $cl \left (  \pi^k_X \circ \Delta^{sm} \circ  \left ( \pi^i_X \times \pi^j_X \right )\big|_{U_b} \right ) =0$ in $H^{*}(U_b, \mathbb Q)$ as on every fiber it gives a Chow-Kunneth decomposition and it is multiplicative cohomologically, by SFP we get     \begin{equation*}
          \pi^k_U \circ \Delta^{sm} \circ  \left ( \pi^i_U \times \pi^j_U \right )  =0\text{ for }k\neq i+j
    \end{equation*}in $A^* (U\times_B U\times_B U)$. This completes the proof. 
\end{proof}
\noindent Now we are ready to prove Theorem \ref{MCK Theorem}. Consider the family $\mathcal C^w_{2,g} \rightarrow \mathcal H^w_g$. This family admits a section given by the marked Weirestrass point. In particular, there is a line bundle $\mathcal L$ on $\mathcal C^w_{2,g}$ of relative degree $1$. The following decomposition is proved in \cite[Theorem 4.3]{hoskins} (see also \cite[Corollary 2.7]{MR1895918}) for any smooth curve $C$ over an \textit{arbitrary} field of $\operatorname{char} 0$ and $L$ an odd degree line bundle. \begin{equation}
    \mathfrak h(U(2,  L)) = \mathfrak h(Sym^{g-1}C) (g-1) \oplus  \bigoplus _{i=0}^{g-2} \mathfrak h(Sym^i C) \otimes \big ( \mathbb Q(i) \oplus \mathbb Q(3g-3-2i) \big ) .
\end{equation} An argument analogous to Theorem \ref{thm:decomp-rel} and the Edidin-Graham approximation gives us the decomposition \begin{equation}\label{MCK motivic decomp}
    \mathfrak h_{\h^w_g} (\mathcal{U}^\h(2,\mathcal L,g))\simeq \mathfrak h _{\h^w_g}( (\mathcal C^w_{2,g})^{(g-1)})(g-1) \oplus  \bigoplus_{i=0}^{g-2} \mathfrak h_{\h^w_g}(\mathcal C^w_{2,g})^{(i)}) \otimes 
 \biggl(\mathbb Q(i) \oplus \mathbb Q(3g-3-2i)  \biggr).
\end{equation}
We are now ready to prove Theorem \ref{MCK Theorem}.
\begin{proof}[Proof of Theorem \ref{MCK Theorem}]
    Since from \cite[Theorem 4.2]{hoskins}, we have $\mathfrak h(U_C(2,d)) \simeq \mathfrak h(Jac(C))\otimes \mathfrak h(U_C(2,L)) $ for any $L\in Pic^d(C)$, where the isomorphism is induced by the pullback of the morphism $U_C(2,L) \times Jac(C) \rightarrow U_C(2,d)$ and hence the motivic isomorphism is an isomorphism of algebra objects. So it suffices to prove the theorem for $U_C(2,\mathcal L)$ for any line bundle $\mathcal L$ with $\deg \mathcal L=1$. 

    \noindent Now Corollary \ref{cor:hyp-weier-powers} shows that $(\mathcal C^w_{2,g})^{\times_{\mathcal H^w_g}n}\rightarrow \mathcal H^w_g$ satisfies SFP for $n\leq 2g+4$. Since $g\leq 7$, $3g-3 \leq 2g+4$. So every relative triple product on the RHS of Equation \ref{MCK motivic decomp} satisfies SFP. Further since $\mathcal C^w_{2,g} \rightarrow \mathcal H^w_g$ admits a relative CK decomposition, so does $(\mathcal C^w_{2,g})^{(k)}\rightarrow \mathcal H^w_g$. Lemma \ref{triple SFP-MCK Lemma} then completes the proof.
\end{proof}
\begin{rem}
    Existence of MCK for moduli spaces of semistable rank 2 odd degree bundles over genus $2$ curves can also be deduced from \cite{MR429897} and Lateveer's \cite[Theorem 4.1]{laterveer-franchetta}. But Theorem \ref{MCK Theorem} for $U_C(2,\mathcal L)$ with $3\leq g(C) \leq 7$ gives new examples of Fano varieties admitting an MCK to the best of the authors' knowledge.  
\end{rem}
\section{Scripts for supporting computations}
\begin{lstlisting}[language=Python]
import sympy as sp
T, g = sp.symbols('T g', positive=True)
X = (1-T**g)/(1-T) + (T**(g+1)-T**(3*g-1))/(1-T**2)
Y = (1-T**(2*g))/(1-T**2) + (T**(2*g-1)-T**(3*g-2))/(1-T)
Z = (1-T**(3*g-sp.Rational(3,2)))/(1-T**sp.Rational(3,2))
Zp = (1+T**(3*g-sp.Rational(3,2)))/(1+T**sp.Rational(3,2))
W = (g*T**(g+1)+1-(g+1)*T**g)/(1-T)**2 \
    + ((g-1)*T**(g+1)+T**(3*g+1)-g*T**(g+3))/(1-T**2)**2
V = (g*T**(2*g+2)+1-(g+1)*T**(2*g))/(1-T**2)**2 \
    + T**(2*g-1)*((g-1)+T**g-g*T)/(1-T)**2
E = X*(1+3*T**(g+2))/(1-T)**4 \
    + Y*(T**(g+4)+3*T**2)/(1-T)**4 \
    + V*T**2/(1-T)**3 \
    - W*T**(g+1)/(1-T)**3 \
    - Z*(T**(g+2)+T**sp.Rational(1,2))/(2*(1-T)*(1-T**sp.Rational(1,2))**3) \
    + Zp*(T**sp.Rational(1,2)-T**(g+2))/(2*(1-T)*(1+T**sp.Rational(1,2))**3)
F = sp.simplify(E*(1-T)**2*(1-T**2)**2*(1-T**3))
F_simplified = sp.simplify(sp.together(F))
sp.factor(F_simplified)
    \end{lstlisting}
\captionof{figure}{ Computations for Hilbert Series of $\mcalu^\h(2,d,g)$}
\label{hilb:code}

    \begin{lstlisting}[language=Python]
ka = E /(1-T)+ T**g*(
    (g+2)*(g+1)/2*T**g
    + sp.summation((g-k+2)*(g-k+1)/2*(T**(g-k)+T**(g+k)), (k, 1, g))
   ) / ((1-T)**2*(1-T**2))
ka_simplified = sp.simplify(sp.together(ka*(1-T)**3*(1-T**2)**2*(1-T**3)))
diff = sp.simplify(ka_simplified - (1-T**(g+1))*(1-T**(g+2))*(1-T**(2*g+3)))
diff
    \end{lstlisting}
   \captionof{figure}{Checking equality of Hilbert Series} 
   \label{conj:code}

\bibliographystyle{abbrv}
\bibliography{ack-h}
\end{document}